\documentclass{article}
\usepackage{latexsym,amsmath,amssymb,amscd,xcolor,bm}
  \newtheorem{notation}{Notation}
    \newtheorem{remark}{Remark}

\begin{document}

\title{An Application of Nash-Moser Theorem to Smooth Solutions of One-Dimensional Compressible Euler Equation with Gravity
}
\author{Cheng-Hsiung Hsu \footnote{Department of Mathematics, National Central University, Chungli 32001, Taiwan. 
          Research supported in part by NSC and NCTS of Taiwan. e-mail:chhsu@math.ncu.edu.tw}, Song-Sun Lin 
\footnote{Department of Applied Mathematics, National Chiao-Tung University, Hsinchu 30050, Taiwan. Research supported in part by Center of S.-T. Yau, NCTU, Taiwan; NSC and NCTS of Taiwan. e-mail:sslin@math.nctu.edu.tw}, Tetu Makino \footnote{Faculty of Engineering, Yamaguchi University, Ube 755-8611, Japan. Research supported in part by Center of S.-T. Yau, NCTU, Taiwan. e-mail:makino@ yamaguchi-u.ac.jp}\\
 and  Chi-Ru Yang \footnote{Mathematics Division, National Center for Theoretical Sciences, Hsinchu 30043, Taiwan. Research supported in part by NCTS of Taiwan. e-mail:cryang@math.cts.nthu.edu.tw}}

\date{\today}

\maketitle

\begin{abstract}
We study one-dimensional motions of polytropic gas governed by the compressible Euler equations. 
The problem on the half space under a constant gravity gives an equilibrium which has free 
boundary touching the vacuum and the linearized approximation at this equilibrium gives time 
periodic solutions. But it is not easy to justify the existence of long-time true solutions for 
which this time periodic solution is the first approximation. The situation is in contrast to 
the problem of free motions without gravity. The reason is that the usual iteration method 
for quasilinear hyperbolic problem cannot be used because of the loss of regularities which causes from
 the  touch with the vacuum. Interestingly, the equation can be transformed to a nonlinear wave equation 
on a higher dimensional space, for which the space dimension, being larger than 4, is related to the 
adiabatic exponent of the original one-dimensional problem. We try to find a family of solutions 
expanded by a small parameter. Applying the Nash-Moser theory, we justify this expansion.
The application of the Nash-Moser theory is necessary for the sake of conquest of 
the trouble with loss of regularities, and the justification of the applicability 
requires a very delicate analysis of the problem.\\


\end{abstract}

\section{Introduction}
\label{intro}
The aim of this paper is to study one-dimensional motions of polytropic gas governed by the compressible Euler equations
\begin{align}\rho_{{t}}+(\rho {u})_{{x}}&=0,\label{Euler1}\\
(\rho{u})_{{t}}+(\rho{u}^2+P)_{{x}}&=- g\rho,\label{Euler2}
\end{align}
for $t,x\ge 0$ subject to the boundary condition
\begin{align}\label{BC}
\rho u|_{x=0}=0.
\end{align}
Here $\rho,u,P$ and $g>0$ are density, velocity, pressure and gravitational acceleration constant respectively.
Equations \eqref{Euler1}$\sim$\eqref{BC} describe the atmosphere on the flat earth $\{x\le 0\}$ moving in one direction under the constant gravitational force downward.
In this work we assume that $P=P(\rho)=A\rho^\gamma$ for some constants $A,\gamma$ such that
 $0<A, 1<\gamma\le 2$. 
Then equilibria of \eqref{Euler1}  and \eqref{Euler2} are of the form
\begin{align}\label{steady}
\bar{\rho}=\left\{
\begin{array}{rl}
A_1(x_+-x)^{\frac{1}{\gamma-1}}, & \mbox{if }0\leq x \leq
x_+,\\
 0,& \mbox{if } x_+<x,
 \end{array}
 \right.
 \end{align}
 where $A_1= ((\gamma -1)g/\gamma
A)^{1/(\gamma-1)}$ and $x_+$ is an arbitrary positive value,
which represents the stratospheric depth.
\medskip

Without loss of generality, we may assume $x_+=1$, $A_1=1$ and $A=1/\gamma$. It can be seen easily
 by scale transformations of the variables.  Since the interface with the vacuum would vary with the time, it
is convenient to transform the equations \eqref{Euler1} and \eqref{Euler2} into the Lagrangian form. More precisely, we introduce the variable $$ m=\int_0^x\rho dx $$ as the
independent variable instead of $x$,  then equations \eqref{Euler1} and \eqref{Euler2} can be transformed   into the following second order equation:
\begin{align}\label{single}
x_{tt}+P_m=-g,
\end{align}
where $P=\gamma^{-1}(x_m)^{-\gamma}.$
Let us fix an equilibrium
\begin{align}\label{steady1}
x=\bar{x}(m)=1-A_2(m_+-m)^{\frac{\gamma-1}{\gamma}},
\quad 0\leq m \leq m_+,
\end{align} where
\begin{align*}
A_2=\Big(\frac{\gamma}{\gamma-1}\Big)^{\frac{\gamma}{\gamma-1}}\ \ \mbox{and}\ \
m_+=\frac{\gamma-1}{\gamma}.
\end{align*}
Then we consider small perturbations
of the equilibrium in \eqref{steady1} by putting $ x(t,m)=\bar{x}(m)+y.$ Under such assumption, the
equation \eqref{single} is reduced to
\begin{align}\label{wave}
y_{tt}-(\gamma\bar{P}G(\frac{1}{\bar{x}_m}y_m))_m=0,
\end{align}where $$
G(v)=\frac{1}{\gamma}(1-(1+v)^{-\gamma}). $$ 
Take $\bar{x}(m)$ as the independent
variable instead of $m$, and write it as $x$.
Then the
equation \eqref{wave} is reduced to
\begin{equation}\label{linear}
y_{tt}-\frac{1}{\bar{\rho}(x)}(\gamma \bar{P}(x)G(y_x))_x=0
\end{equation}
for $0<x<1$ and the boundary condition is
\begin{equation}\label{BC1}
y|_{x=0}=0.
\end{equation}
Here $$
\bar{\rho}(x)=(1-x)^{\frac{1}{\gamma-1}}
\quad\mbox{and}\quad
\bar{P}(x)=\frac{1}{\gamma}(1-x)^{\frac{\gamma}{\gamma-1}}. $$
Equation \eqref{linear} is an apparently quasilinear hyperbolic equation. But it 
has a singularity at $x=1$. Due to the singularity, the investigation for the existence of
 time periodic solutions becomes a difficult and challenging problem. To the best of our knowledge, 
the existence problem of time periodic solutions is still open. 

\medskip
For the sake of comparison, let us recall the results of \cite{HLM2},
which considered the following simplified quasilinear wave equation
\begin{align}\label{simplify}
\left\{
\begin{array}{ll}
y_{tt}-(G(y_x))_x=0\ \ \mbox{for }0<x<1,\medskip\\
y(t,0)=y(t,1)=0.
\end{array}
\right.
\end{align}
This problem is derived from the Euler equations
$$\rho_t+(\rho u)_x=0, $$
$$(\rho u)_t+(\rho u^2+P)_x=0, $$
and the boundary condition
$$\rho u|_{x=0}=\rho u|_{x=L}=0, $$
for which the equilibria are positive constant densities. Hence there are no troubles
caused by contact with vacuum.
For any fixed arbitrarily long time, \cite{HLM2} shows that there are smooth small amplitude solutions of 
the problem (10) for which the periodic solutions of the linearized equation
 are the first-order approximation. 
This result was established by the usual iteration method for quasi-linear wave equations.

Therefore, similarly, we want to find smooth solutions for which a time periodic solution of the linearized
equation around an equilibrium is the first approximation even for the present
problem (8)(9).   
However, contrary to the case without gravity, the usual iteration method 
for quasilinear hyperbolic problem cannot be applied directly to 
the present problem because of the loss of regularities which causes from the touch with the vacuum. 
In this work we shall apply Nash-Moser theorem to establish long time existence 
of smooth solutions near time-periodic solution of the linearized equation.

More precisely speaking, we introduce the variable
\begin{equation}
z=1-x
\end{equation}
and small parameter $\varepsilon$, and we shall construct
approximate solutions of the form
$$\sum_{k=1}^K y_k(t,z)\varepsilon^k,$$
where $y_k(t,z)$ are entire functions of $t$ and $z$, while $y_1(t,z)$ is a non-trivial time
periodic solution of the linearized equation.

Then our aim is to find a true smooth solution $y(t,z)$ of (8)(9)
on $0\leq t\leq T$ and $0\leq z\leq 1$, for arbitrarily fixed $T$, such that 
$$y(t,z)=\sum_{k=1}^Ky_k(t,z)\varepsilon^k +O(\varepsilon^{K+1}).$$
Of course for large $T$ we should restrict $\varepsilon$ sufficiently small.
Then 
$$x(t,m)=\bar{x}(m)+y(t,1-\bar{x}(m)) $$
is a solution in the Lagrangian variable and the corresponding density distribution $\rho=\rho(t,x)$,
where $x$ denotes the original Euler coordinate, satisfies
$$\rho(t,x)>0 \qquad \mbox{for}
\qquad 0\leq x<x_F(t)$$
and
$$\rho(t,x)=0\qquad \mbox{for}\qquad x_F(t)\leq x,$$
where
$$x_F(t)=1+y(t,0)$$
is the position of the free boundary. Since $y(t,z)$ is smooth on $0\leq z\leq 1$,
we have
$$\rho(t,x)=C(t)(x_F(t)-x)^{\frac{1}{\gamma-1}}(1+O(x_F(t)-x)), \qquad (x<x_F(t))$$
and 
$$\frac{\partial}{\partial x}\Big(\frac{dP}{d\rho}\Big)=\frac{\partial}{\partial x}\rho^{\gamma-1}
\rightarrow -C(t) $$
at $x \rightarrow x_F(t)-0$. This condition is that of ``physical vacuum boundary" so called by
the most recent works 
\cite{JM}(2009) and \cite{CS}(2011). This concept can be traced back to \cite{Liu}(1996),
\cite{LiuY}(2000), and \cite{Y}(2006). Hence 
we can say that our purpose is to find long-time smooth solutions with 
``physical vacuum boundary". But \cite{JM} and \cite{CS} are interested in short-time solutions 
to
the initial value problem for
 the case without external force. So the motivation, methods and results are different from those of this
work.

\medskip
Now we have introduced the variable
$$z=1-x.$$ 
Moreover it is convenient to introduce the parameter
\begin{equation}
\gamma=1+\frac{2}{N-2}.
\end{equation}
Then the assumption $1<\gamma\le2$ is equivalent to that  $4\le N<\infty$. 
Hence, we assume $N\ge 4$ in the following of this wok.
Moreover, the equation \eqref{linear} turns out to be
\begin{equation}\label{laplace}
\frac{\partial^2y}{\partial t^2}
-\triangle y=
G_I(v)\triangle y +G_{II}(v),
\end{equation}
where
\begin{align}
\triangle &:= z\frac{\partial^2}{\partial z^2}+\frac{N}{2}\frac{\partial}{\partial z},
\quad
v=-\frac{\partial y}{\partial z},\\
G_I(v)&:=DG(v)-1=-\frac{2N-2}{N-2}v+[v]_2,\\
G_{II}(v)&:=\frac{N}{2}(vDG(v)-G(v))=-\frac{N(N-1)}{2(N-2)}v^2+[v]_3
\end{align}
and
$[v]_q$ denotes a convergent power series of the
form $\sum_{j\geq q}a_jv^j$.

 If we introduce the variable $r$ by
\[z=1-x=\frac{r^2}{4}\]
then
\[\triangle= z\frac{\partial^2}{\partial z^2}+\frac{N}{2}\frac{\partial}{\partial z}=\frac{\partial^2}{\partial r^2}+\frac{N-1}{r}\frac{\partial}{\partial r}
\]
is the radial part of the Laplacian operator on the $N$-dimensional Euclidean space ${\mathbb R}^N$,
provided that $N$ is an integer. But we shall not assume that $N$ is an integer
in this work.

\medskip
Here we would like to spend few words to explain why
the usual iteration does not work although the equation (13) is apparently quasi-linear.
For the sake of simplicity, let us assume $N$ is an integer. Then
a smooth function $y$ of $z$ can be regarded as a smooth function of
$r=\|\vec{x}\|= (\sum_j (x_j)^2)^{1/2}$, where $\vec{x} \in \mathbb{R}^N$.
Since $y$ is smooth and spherically symmetric, we can assume that
$\partial y/\partial r=0$ at $r=0$ and
$$-v=\frac{\partial y}{\partial z}=\frac{2}{r}\frac{\partial y}{\partial r} \rightarrow
2\frac{\partial^2y}{\partial r^2}\Big|_{r=0}$$
as $r \rightarrow 0$. In other words, $v=-\partial y/\partial z$ is not of the 
first order, but of the second order, which is of
the same order as the principal part $\displaystyle \triangle y
=z\frac{\partial^2y}{\partial z^2}+\frac{N}{2}\frac{\partial y}{\partial z}$. So,
the loss of
regularities cannot be recovered
by one step of solving a (linear) wave equation. This is the reason why
we try an application of the Nash-Moser theory. Note that
this trouble comes from $z=0$, that is, from the touch with vacuum at the free boundary.

\section{Preparatory analysis of linear problems}
First let us consider the linearized problem of (13):
\begin{equation}
 y_{tt}-\triangle y=0, \qquad y|_{z=1}=0. 
\end{equation}
In \cite{HLM1} we showed that (17) admits a time periodic solution
\begin{equation}\label{y1}
 y=y_1=\sin(\sqrt{\lambda_{n}}t+\theta)\Phi_{\frac{N}{2}-1}(\lambda_{n}z),
\end{equation}
where $\theta$ is a constant, $\lambda_{n}$ is the eigenvalues of the operator $-\triangle $
with the Dirichlet boundary condition, and
$$\Phi_{\frac{N}{2}-1}(X)=
\sum_{k=0}^{\infty}\frac{(-1)^k}{k!\Gamma(\frac{N}{2}-1+k+1)}X^k $$
is an entire function such that $\Phi_{\frac{N}{2}-1}(\lambda_n)=0$.
In fact,
$$\lambda_n=\frac{1}{4}{(j_{\frac{N}{2}-1,n})^2},$$
where $j_{\frac{N}{2}-1,n}$ is the $n$-th positive zero of
the Bessel function $J_{\frac{N}{2}-1}$, and
$$J_{\nu}(\zeta)=\Big(\frac{\zeta}{2}\Big)^{\nu}
\Phi_{\nu}\Big(\frac{\zeta^2}{4}\Big).$$

More precisely speaking, 
we consider the Hilbert space $\mathfrak{X}$ which consists of
functions of $0\le z\le1$ endowed with the inner
product
$$
(y_1|y_2)_{\mathfrak{X}}:=\int_0^1
y_1(z)\overline{y_2(z)}z^{\frac{N}{2}-1}dz.
$$
The self-adjoint operator $T=-\triangle$ with boundary condition is defined on
\begin{align*}\mathfrak{D}(T)=
\{ y \in \mathfrak{X} \  | &\  \exists\ \eta_n \in C_0^{\infty}(0,1)\  \mbox{such that}\
\eta_n\rightarrow y\ \mbox{in}\ \mathfrak{X},\
 Q[\eta_n-\eta_m]\rightarrow 0 \\ &\ \mbox{as }m,n\rightarrow\infty,\ \mbox{and} \ -\triangle y \in \mathfrak{X} \ \mbox{in distribution sense} \}.
\end{align*}
Here
$$Q[\eta]:=\int_0^1\Big|\frac{d\eta}{dz}\Big|^2z^{\frac{N}{2}}dz$$
and
``$-\triangle y =f \in \mathfrak{X}$ in distribution sense" means that
for any $\eta \in C_0^{\infty}(0,1)$ there holds
$$
(y|-\triangle \eta )_{\mathfrak{X}}=(f|\eta)_{\mathfrak{X}}.
$$
By \cite{HLM1}, we have
$$\mathfrak{D}(T)=\{ y\in C(0,1]
\ | \ y \in \mathfrak{X}, \ y(1)=0\ \mbox{and}\ -\triangle y \in \mathfrak{X}
\ \mbox{in distribution sense}\}
$$
and the spectrum of $T$ consists of simple eigenvalues $\lambda_1<\lambda_2<\cdots$,where 
$\lambda_n=j^2_{\nu,n}/4$.

\medskip
Moreover, we consider the problem
\begin{align}\label{elliptic}
-\lambda y-\triangle y&=f(z), \quad z\in(0,1).
\end{align}
Here $\lambda \geq 0$ and $f$ are given.

\newtheorem{proposition}{Proposition}
\begin{proposition}
The inverse $T^{-1}$ is a compact operator.
\end{proposition}
Proof.
If $f\in\mathfrak{X}$, the solution of the problem \eqref{elliptic} with $\lambda=0$ is given by the formula
\begin{eqnarray*}
y(z)=\frac{2}{N-2}\Big(
\int_z^1f(\zeta)d\zeta +
z^{-\frac{N}{2}-1}\int_0^z
f(\zeta)\zeta^{\frac{N}{2}-1}d\zeta -\int_0^1f(\zeta)\zeta^{\frac{N}{2}-1}d\zeta \Big).
\end{eqnarray*}
Since
\begin{align*}
\int_z^1|f(\zeta)|d\zeta &\leq \sqrt{\int_z^1\zeta^{-\frac{N}{2}-1}d\zeta}
\sqrt{\int_z^1|f(\zeta)|^2\zeta^{\frac{N}{2}-1}d\zeta}\\
&\leq \left\{\begin{array}{rl}\sqrt{\frac{2}{N-4}}z^{-\frac{N-4}{4}}||f||_{\mathfrak{X}},&\mbox{if } N>4,\medskip \\
|\log z|||f||_{\mathfrak{X}}, &\mbox{if } N=4
\end{array}
\right.
\end{align*}
and
\begin{align*}
\int_0^z|f(\zeta)|\zeta^{\frac{N}{2}-1}d\zeta\leq \sqrt{\frac{2}{N}}z^{\frac{N}{4}}||f||_{\mathfrak{X}},
\end{align*}
we see that
\begin{eqnarray*}
|y(z)|\leq \left\{\begin{array}{rl}
C z^{-\frac{N-4}{4}}||f||_{\mathfrak{X}},&\mbox{if } N>4, \medskip\\
C |\log z|\cdot\|f\|_{\mathfrak{X}},&\mbox{if }\ N=4
\end{array}
\right.
\end{eqnarray*}
for some constant $C>0$. Moreover, we have
$$\frac{dy}{dz}=-z^{-\frac{N}{2}}\int_0^zf(\zeta)\zeta^{\frac{N}{2}-1}d\zeta$$
and which implies
$$\Big|\frac{dy}{dz}  \Big|\leq \sqrt{\frac{2}{N}}z^{-\frac{N}{4}}||f||_{\mathfrak{X}}. $$
Therefore, Ascoli-Arzel\`{a}'s theorem implies that a sequence $y_n$
converges on each compact subset of $(0,1]$ when
$f$ is confined in a bounded set of $\mathfrak{X}$.
On the other hand, since
\begin{align*}
\int_0^{\delta}|y(z)|^2z^{\frac{N}{2}-1}dz
\leq \left\{
\begin{array}{rl}
C {\delta}^2||f||_{\mathfrak{X}}^2, &\mbox{if } N>4, \medskip\\
C {\delta}^2|\log \delta|^2\cdot\|f\|^2_{\mathfrak{X}},&\mbox{if }\ N=4,
\end{array}
\right.
\end{align*}
we see that $y_n$ converges in $\mathfrak{X}$, too. $\square$

Therefore, $T$ is a self-adjoint operator whose inverse is compact and the following assertion holds.
See, eg., \cite{DS}.

\begin{proposition}\label{ttt}
If $\lambda\geq 0$, then the range $\mathcal{R}(-\lambda+T)$ is closed and
$$\mathcal{R}(-\lambda+T)=\mathcal{N}(-\lambda +T)^{\perp}.$$
Thus, if $\lambda=\lambda_n$ is an eigenvalue with an eigenfunction $\phi_n$,
then the problem of \eqref{elliptic}
admits a solution $y$ in $\mathfrak{X}$ if and only if
$$(f|\phi_n)_{\mathfrak{X}}=0.$$
\end{proposition}

\begin{proposition}\label{prop1}
If $f(z)$ is an entire function, then there is an entire function $y(z)$ which
solves the equations of \eqref{elliptic}.
\end{proposition}
Proof.
Let
\begin{align}
f(z)=\sum_{k=0}^{\infty} c_kz^k.
\end{align}
Let $R$ be an arbitrarily large positive number. Since $f$ is an entire
function, there is a constant $M$ such that
$|c_k|\leq M/R^k$ for all $k$.
We seek a solution $y(z)$ of \eqref{elliptic} in the form
\begin{equation}\label{entire}
y=\sum_{k=0}^{\infty}a_kz^k.
\end{equation}
Substituting \eqref{entire} into \eqref{elliptic} and comparing the coefficients,
we have the formula
$$a_{k+1}=-\frac{\lambda a_k+c_k}{(k+1)(k+\frac{N}{2})}.$$
Taking $a_0$ arbitrarily, we claim that there is a constant $\bar{M}>0$ such that $|a_k|\leq
\bar{M}/R^k$ for all $k$. Suppose $k>\max\{R,\lambda+1\}$ and
$|a_k|\leq {M}^\prime/R^k$, then
$$|a_{k+1}|\leq \frac{\lambda|a_k|+|c_k|}{k^2}
\leq\frac{(\lambda+1)(M+{M}^\prime)}{kR^{k+1}}\le \frac{M^\prime}{R^{k+1}}$$
provided that
$$\frac{(\lambda+1)(M+ M^\prime)}{k} \leq  M^\prime.$$
Hence the claim follows and
the radius of convergence of $\sum a_kz^k$ is larger than $R$.  $\square$

\begin{proposition}\label{prop34}
Suppose $\lambda >0$ and $f$ is an entire function, then
any solution of \eqref{elliptic} in $\mathfrak{X}$ is an entire function.
\end{proposition}
Proof.
 The homogeneous equation
$-\lambda y-\triangle y=0$
admits a pair of linearly independent
solutions
$y_1(z)=\Phi_{\frac{N}{2}-1}(\lambda z)$ and
$y_2(z)$ such that
$$y_2(z) \sim (\lambda z)^{-\frac{N}{2}-1} \ \mbox{as } z\rightarrow 0.$$
In fact, if we take the change of variables
$$\lambda z=\frac{r^2}{4} \ \ \mbox { and }\ \ y=r^{-\nu}w, $$
then the equation $-\lambda y-\triangle y=0$ turns out to be
the following Bessel equation:
$$\frac{d^2w}{dr^2}+\frac{1}{r}\frac{dw}{dr}+
\Big(1-\frac{\nu^2}{r^2}\Big)w=0. $$
If $\nu$ is not an integer then $J_{\nu}$ and $J_{-\nu}$ are linearly independent solutions. On the other hand, if
$\nu(\not =0)$ is an integer, then $J_{\nu}$ and
the Bessel function of the second kind $Y_{\nu}$ of the form
\begin{align*}
Y_{\nu}(r)=
\frac{2}{\pi}J_{\nu}(r)\log \frac{r}{2}  -&\frac{1}{\pi}\Big(\frac{r}{2}\Big)^{\nu}\sum_{k=0}^{\infty}
\frac{(-1)^k(\Psi(k+1)+\Psi(\nu+k+1))}{\nu !(\nu +k)!}\Big(\frac{r^2}{4} \Big)^k \\
-&\frac{1}{\pi}\Big(\frac{r}{2}\Big)^{-\nu}\sum_{k=0}^{\nu-1}
\frac{(\nu-1-k)!}{k!}\Big(\frac{r^2}{4} \Big)^k
\end{align*}
are linearly independent solutions.
Here $\Psi(x):=D\Gamma(x)/\Gamma(x)$. See \cite{W}.
Since $N\geq 4$, we see $y_2$ does not belong to $\mathfrak{X}$.
On the other hand, there is an entire function $y=\psi_0(z)$ which satisfies \eqref{elliptic} due to Proposition \ref{prop1}.
Of course $\psi_0 \in \mathfrak{X}.$ Thus any solution $y(z)$ of \eqref{elliptic} can be written as
$$y(z)=\psi_0(z)+C_1y_1(z)+C_2y_2(z),$$
in which, if $y(z) \in \mathfrak{X}$, then $C_2=0$, and therefore, $y(z)$ is an entire function, too. $\square$

\section{Formal solution expanded as power series of parameters}
Now we construct formal power series solution of \eqref{laplace}. Let us fix
a non-trivial solution
$$y_1=\sin(\sqrt{\lambda_{n_0}}t+\theta_0)
\phi_{n_0}(z)$$
of the linearized problem,
where
\begin{align}\label{normal}
\phi_n(z)=\frac{\Phi_{\frac{N}{2}-1}(\lambda_nz)}{||\Phi_{\frac{N}{2}-1}(\lambda_nz)||_{\mathfrak{X}}}
\end{align}
is the normalized eigenfunction in the Hilbert space $\mathfrak{X}$.
According to the result of \cite{HLM1}, we know that  $(\phi_n)_{n=1,2,\cdots}$ forms a 
complete orthonormal system in $\mathfrak{X}$. Note that $\phi_n(z)$ is an entire function of $z$. Our 
purpose is to find a formal solution of \eqref{laplace} of
the form
\begin{equation}\label{taylor}
y(t,z)=\sum_{k=1}^{\infty}y_k(t,z)\varepsilon^{k},
\end{equation}
where $\varepsilon$ stands for a small parameter. Substituting \eqref{taylor} 
into the equation \eqref{laplace} and comparing the coefficients,
we get the following sequence of linear equations
\begin{align}\label{seq}
\Big(\frac{\partial^2}{\partial t^2}-\triangle\Big)y_k=&
\sum_{1 \leq \ell, \  j_1+\cdots+j_\ell+j=k} G_{I\ell}v_{j_1}\cdots v_{j_\ell}\triangle y_j +\nonumber\\
&+\sum_{2 \leq \ell, \  j_1+\cdots+j_\ell=k} G_{II\ell}v_{j_1}\cdots v_{j_\ell} ,
\end{align}
where
$$v_j=-\frac{\partial y_j}{\partial z},\ \ G_{I}(v)=\sum_{1\leq \ell}G_{I\ell}v^\ell\ \ \mbox{and}\ \
G_{II}(v)=\sum_{2\leq \ell}G_{II\ell}v^\ell.
$$
Starting from the fixed $y_1$, we can solve the equations \eqref{seq} with the boundary condition
$y_k(1)=0$ successively.
\subsection{Solution for $k=2$}

The equation of \eqref{seq} for $k=2$ is in the form
\begin{align}\label{seq2}
\Big(\frac{\partial^2}{\partial t^2}-\triangle\Big)y_2=G_{I1}v_1\triangle y_1+G_{II2}(v_1)^2=-\frac{2(N-1)}{N-2}(\triangle y_1+\frac{N}{4}v_1)v_1.\quad
\end{align}
Since $y_1$ is an entire function,  we can write the right-hand side
of \eqref{seq2} by the form
$$f_0(z)+(\cos2\Theta)f_1(z), $$
where $\Theta :=\sqrt{\lambda_{n_0}}t+\theta_0,$
 $f_0$ and $f_1$ are entire functions of $z$.
 Here we have used
$$\sin^2\Theta =\frac{1}{2}-\frac{1}{2}\cos 2\Theta.$$
Then solutions of the problem
\begin{align}\label{seq2a}
-\triangle w= f_0(z),\quad w|_{z=1}=0
\end{align}
can be represented by
\begin{align}\label{solution}
w=-\frac{2}{N-2}\int_0^z\Big(1-\Big(\frac{\zeta}{z}\Big)^{\frac{N}{2}-1}\Big)f_0(\zeta)d\zeta
+\frac{2}{N-2}\int_0^1(1-\zeta^{\frac{N}{2}-1})f_0(\zeta)d\zeta.
\end{align}
On the other hand, we consider the problem
\begin{align}\label{seq2b}
\Big(\frac{\partial^2}{\partial t^2}-\triangle\Big)w=(\cos2\Theta)f_1(z), \quad w|_{z=1}=0.
\end{align}
We need to consider the following two cases:\medskip

\qquad Case-1: $4\lambda_{n_0}$ is not
an eigenvalue; \medskip

\qquad Case-2: there is an eigenvalue $\lambda_q=4\lambda_{n_0}$.\medskip

First, we consider the Case-1. 
 Then \eqref{seq2b} has a solution of the form
$w(t,z)=(\cos2\Theta)W(z),$
where $W(z)$ satisfies
\begin{align}\label{seq2c}
\big(-4\lambda_{n_0}-\triangle \big)W=f_1(z), \quad W|_{z=1}=0.
\end{align}
According to Proposition \ref{prop1}, the first equation of \eqref{seq2c} has a solution $W_0(z)$, which is an entire function of $z$ such that
$W_0(0)=1$. Then,  for any constant $C$,
$$W(z)=W_0(z)+C\Phi_{\frac{N}{2}-1}(4\lambda_{n_0}z)$$
is a solution of \eqref{seq2c}, too. Since $4\lambda_{n_0}$ is not an eigenvalue,
we have $\Phi_{\frac{N}{2}-1}(4\lambda_{n_0})\not= 0$. Therefore, we can choose
$C$ so that
\begin{center}
$W(1)=W_0(1)+C\Phi_{\frac{N}{2}-1}(4\lambda_{n_0})=0$,
\end{center}
 i.e. $W(z)$ satisfies the boundary value condition.\medskip

Next, we consider Case-2, i.e. $\lambda_q=4\lambda_{n_0}$ for some integer $q$. We guess that this case
could not happen actually for $N\geq 4$. 
More generally, we have \medskip 

\noindent{\textbf{Conjecture}}{\textit{
Let $\nu\geq 1$ and $\theta$ be a positive zero of the Bessel function
$J_{\nu}$. Then $J_{\nu}(L\theta)\not=0$ for any integer $L\geq 2$.}}\medskip

\noindent (Note that the conclusion is not the case if $\nu=1/2$, for which
$J_{1/2}(r)=\displaystyle \sqrt{\frac{2}{\pi r}}\sin r$.) However 
we have not yet verified this conjecture.  Therefore we should consider Case-2. By Proposition 2,
there is a solution $W_1(z)$ of
$$(-\lambda_q-\triangle )W_1=\tilde{f}_1(z):=
f_1(z)-(f_1|\phi_q)_{\mathfrak{X}}\phi_q(z), $$
which is entire and satisfies the boundary condition.
Then it is easy to see that
$$w=(\cos 2\Theta)W_1-
\frac{1}{2\sqrt{\lambda_q}}t\cdot\sin 2\Theta \cdot(f_1|\phi_q)_{\mathfrak{X}}\phi_q(z) $$
satisfies \eqref{seq2b}.\medskip

Summing up, we have a solution $y_2$ of the form
\begin{equation}\label{2nd}
y_2(t,z)=y_{20}(z)+(\cos 2\Theta)y_{21}(z),
\end{equation}
for Case-1, or
$$y_2(t,z)=y_{20}(z) +(\cos2\Theta )y_{21}(z) +
t(\sin 2\Theta)y_{22}(z), $$
for Case-2,
where $y_{20}, y_{21}$ and $y_{22}$ are entire functions of $z$. 

\medskip
Suppose Case-1. Then the 2nd order approximate solution
$$y^{(2)}(t,z)=\varepsilon y_1(t,z)+\varepsilon^2y_2(t,z)$$
is a time-periodic solution with period $\Omega=2\pi/\sqrt{\lambda_{n_0}}$.

According to {\rm \cite{HLM1}} ,
we know that any non-trivial true time-periodic
solution of \eqref{simplify}, if exists, with period $T$ should satisfy
\begin{align}\label{necessary}
\frac{1}{T}\int_0^{T}y(t,x=1)dt >0
\end{align}
as an effect of nonlinearity.
This is true for the approximate solution $y^{(2)}$, since
$$\frac{1}{T}\int_0^{T}y^{(2)}(t,x=1)dt=\varepsilon^2 y_{20}(0)$$
where
$$y_{20}(0)=\frac{N-1}{2(N-2)}\int_0^1
\Big(\frac{d\phi_{n_0}}{dz} \Big)^2dz >0.$$
Using integration by parts, the form of $y_{20}(0)$ can be checked by a tedious but direct computation.

We do not know whether Case-2 actually happens and the resonance could occur or not. This
is an interesting open problem.
\subsection{Solution for  $k=3$}

For the sake of simplicity, we assume Case-1 for $k=2$.
Fixing $y_2$ of the form \eqref{2nd}, we see that the right-hand side of the
equation $\eqref{seq}$ is of the form
$$(\sin\Theta)g_1(z)+(\sin 3\Theta)g_3(z), $$
where
$g_1$ and $g_3$ are entire functions of $z$. Here
we have used
$$\sin^3\Theta = \frac{3}{4} \sin\Theta -\frac{1}{4}\sin 3\Theta. $$
The equation
\begin{align}\label{seq3}
\Big(\frac{\partial^2}{\partial t^2} -\triangle  \Big)w=(\sin\Theta)g_1(z)
\end{align}
has a solution of the form
$$w_1(t,z)=(\sin\Theta)W_1(z)-
\frac{1}{2\sqrt{\lambda_{n_0}}}t\cdot(\cos\Theta)\cdot (g_1|\phi_{n_0})\phi_{n_0}(z), $$
where $W_1(z)$ is an entire function which solves the following equation
\begin{align}
\Big(-\lambda_{n_0}-\triangle  \Big)W_1=\tilde{g}_1:=g_1-(g_1|\phi_{n_0})\phi_{n_0}.
\end{align}
Then we can claim that $w_1(t,z)$ satisfies the boundary condition for any $W_1(z)$
with arbitrary $W_1(0)$. In fact, by Proposition \ref{ttt}, there is a solution $W$ in $\mathfrak{X}$ satisfying the equation
\begin{align*}
(-\lambda_{n_0}-\Delta)W=\tilde g_1
\end{align*}
and the boundary condition.  Then $U:=W_1-W$ belongs to $\mathfrak{X}$ and satisfies the homogeneous equation
\begin{align*}
(-\lambda_{n_0}-\Delta)U=0.
\end{align*}
As in the proof of Proposition \ref{prop34}, there is a constant $C$ such that $U=C\phi_{n_0}$. Hence, 
$W_1=W+U=W+C\phi_{n_0}$ satisfies the boundary condition.\medskip

On the other hand, we consider the problem
\begin{align}\label{513}
\Big(\frac{\partial^2}{\partial t^2}-\triangle\Big)w=(\sin 3\Theta)g_3(z), \quad w|_{z=1}=0.
\end{align}
Similar to the discussion of $k=2$, we need to consider the following two cases:\medskip

\qquad Case-1.1: $9\lambda_{n_0}$ is not an eigenvalue; \medskip

\qquad Case-1.2: there is an eigenvalue $\lambda_q=9\lambda_{n_0}$.\medskip

 Let us consider the Case-1.1 for simplicity.
The problem \eqref{513} has a solution of the form
$w_3(t,z)=(\sin3\Theta)W_3(z),$
where $W_3(z)$ satisfies
\begin{align}\label{seq3c}
\big(-9\lambda_{n_0}-\triangle \big)W_3=g_3(z), \quad W_3|_{z=1}=0.
\end{align}
According to Proposition \ref{prop1}, the first equation of \eqref{seq3c} has a solution $\bar{W}_3(z)$, which is an entire function of $z$ such that
$\bar{W}_3(0)=1$. Then, there exists a constant $C\ne 0$ such that
$$W_3(z)=\bar{W}_3(z)+C\Phi_{\frac{N}{2}-1}(9\lambda_{n_0}z)$$
is a solution of \eqref{seq3c}, too. Since $9\lambda_{n_0}$ is not an eigenvalue, we have
$\Phi_{\frac{N}{2}-1}(9\lambda_{n_0})\not= 0$. Therefore, we can choose
$C$ so that
\begin{center}
$W_3(1)=\bar{W}_3(1)+C\Phi_{\frac{N}{2}-1}(9\lambda_{n_0})=0$,
\end{center}
 i.e. $W_3(z)$ satisfies the boundary value condition.
 Thus,  we have a solution $y_3(t,z)$ of \eqref{seq} of the form
\begin{equation}\label{3rd}
y_3(t,z)=
C_3t(\cos\Theta)\phi_{n_0}(z)+
(\sin\Theta)y_{31}(z)+(\sin3\Theta)y_{33}(z),
\end{equation}
where $y_{31}$ and $y_{33}$ are entire functions of $z$ and $C_3$ is a constant.\medskip

Now we can ask whether the approximate solution $y_3$ given by \eqref{3rd}
is time-periodic or not, or in other words, whether
the resonance occurs actually or not. It depends
on whether
$$C_3=-\frac{1}{2\sqrt{\lambda_{n}}}(g_1|\phi_{n_0})$$
vanishes or not.
We guess that $C_3\not= 0$, but we have not yet verified it.
\subsection{Solutions for $k>3$}

Similar to the computations for the conclusions of
previous subsections, we can determine solutions
$y_k$ of \eqref{seq} for $k>3$ successively in the form
\begin{equation}
y_k(t,z)=
\sum_{M\leq k-1,\ L\leq k}t^M((\cos L\Theta)V_{k,L,M}(z)+(\sin L\Theta)W_{k,L,M}(z)),
\end{equation}
where $V_{k,L,M}$ and $W_{k,L,M,}$ are entire functions of $z$. In order to prove it, 
we need the following
lemma.
\newtheorem{lemma}{Lemma}
\begin{lemma}
 If $f(z)$ is an entire function of $z$,
then the problem
\begin{align}
\Big(\frac{\partial^2}{\partial t^2}-\triangle   \Big)y=t^M(\cos L\Theta) f(z)\ \
\mbox{or}\ \
\Big(\frac{\partial^2}{\partial t^2}-\triangle   \Big)y=t^M(\sin L\Theta) f(z)
\end{align}
admits a solution of the form
\begin{align} y(t,z)=\sum_{m=0}^{M+1}t^m(\cos L\Theta)Y_m(z)\ \ \mbox{or}\ \
y(t,z)=\sum_{m=0}^{M+1}t^m(\sin L\Theta)Y_m(z),
\end{align}
respectively. Here $Y_m(z)$ are entire functions of $z$.
\end{lemma}
Proof.
The proof can be done easily when  $\sqrt{\lambda_m}\ne L\sqrt{\lambda_{n}}$ for all $m,n\in{\mathbb Z}^+$. Here we only consider the case that there is a $q\in{\mathbb Z}^+$ such that $\sqrt{\lambda_q}=L\sqrt{\lambda_{n}}$.
This happens at least if $L=1$. In this case a solution $\psi(z)$ of
$$(-\lambda_q-\triangle)\psi =\tilde{f}:=f-(f|\phi_q)\phi_q$$
which is an entire function of $z$ satisfies the boundary condition.
Then the problem
$$\Big(\frac{\partial^2}{\partial t^2}-\triangle  \Big)y=t^M(\cos L\Theta)f+(-2ML\sqrt{\lambda}t^{M-1}\sin L\Theta +M(M-1)t^{M-2}\cos L\Theta )\psi$$
admits a solution of the form
$$y=t^M(\cos L\Theta)\psi +
A(t)(f|\phi_q)\phi_q(z),$$
where
$$A(t)=\frac{1}{\sqrt{\lambda_q}}
\int_0^t\sin\sqrt{\lambda_q}(t-\tau)\tau^M\cos L\Theta(\tau)d\tau,
$$
in which
$$L\Theta(\tau)=L(\sqrt{\lambda_{n_0}}\tau +\theta_0)=\sqrt{\lambda_q}\tau+L\theta_0.$$
In fact, $A(t)$ is a solution of the equation
$$\frac{d^2A}{dt^2}+\lambda_q A=t^M\cos L\Theta.$$
We see
$$A(t)=\frac{1}{2\sqrt{\lambda_q}}(\sin L\Theta)\frac{t^{M+1}}{M+1}+O(t^M).$$
Then, using the
mathematical induction with respect to $M$, the assertion of the Lemma follows . $\square$

\section{Existence of smooth solutions}

In this section, we will prove the existence of smooth solutions of \eqref{laplace},
using the Nash-Moser theorem.
\medskip
In Section 3, we constructed the approximate solutions
\begin{align}\label{map}
y^{(K)}(t,z):=\sum_{k=1}^{K}y_k(t,z)\varepsilon^k.
\end{align}
Fixing an arbitrarily large $T$ and an integer $K$, we want to find a solution $y(t,z)$ of
the original problem of the form
$$y(t,z)=y^{(K)}+\varepsilon^Kw(t,z)$$
on the time interval $0\leq t\leq T$. First, we derive the problem of $w$. By \eqref{seq}, it follows that $y^{(K)}$ satisfies
\begin{align*}
\frac{\partial^2 y^{(K)}}{\partial t^2}-\triangle y^{(K)}=&\sum_{k=1}^K(\sum_{j_1+\cdots+j_\ell+j=k}G_{I\ell}v_{j_1}\cdots v_{j_\ell}\triangle
y_j)\varepsilon^k+\\
&+\sum_{k=1}^K(\sum_{j_1+\cdots+j_\ell=k}G_{II\ell}v_{j_1}\cdots v_{j_\ell})\varepsilon^k.
\end{align*}
Let us denote $v^{(K)}:=-\partial  y^{(K)}/\partial z$\ \  and\ \  ${P}:=-{\partial w}/{\partial z}$.
Then $w$ satisfies
\begin{align}\label{weq}
\frac{\partial^2 w}{\partial t^2}-\triangle w=G_{I}(v^{(K)}+\varepsilon^{K}P)\triangle w+F_I+F_{II},
\end{align}
where
\begin{align}
\varepsilon^{K+1}F_I=&G_I(v^{(K)}+\varepsilon^KP)\triangle y^{(K)}
-\sum_{k=1}^K\Big(\sum_{j_1+\cdots+j_{\ell}+j=k}
G_{I \ell}v_{j_1}\cdots v_{j_{\ell}}\triangle y_j\Big)\varepsilon^k, \label{g1}\\
\varepsilon^{K+1}F_{II}
=&G_{II}(v^{(K)}
+\varepsilon^KP) -\sum_{k=1}^K\Big(\sum_{j_1+\cdots j_{\ell}=k}
G_{II\ell}v_{j_1}\cdots v_{j_{\ell}}
\Big)\varepsilon^k.\label{g2}
\end{align}
Let us denote
\begin{align}
\varepsilon a(t,z,P,\varepsilon):=& G_I(v^{(K)}(t,z)+\varepsilon^K P), \label{aa}\\
b(t,z,P,\varepsilon):=& -(F_I+F_{II})+(F_I+F_{II})|_{P=0},\label{bb} \\
c(t,z,\varepsilon):=& (F_I+F_{II})|_{P=0}.\label{cc}
\end{align}
Then equation \eqref{weq}  can be written as
\begin{equation}\label{yyy}
\frac{\partial^2 w}{\partial t^2}-\Big(1+\varepsilon a(t,z,-\frac{\partial w}{\partial z},\varepsilon)\Big)
\triangle w+
\varepsilon b(t,z,-\frac{\partial w}{\partial z}, \varepsilon)=\varepsilon c(t,z,\varepsilon).
\end{equation}
Note that $a(t,z, P, \varepsilon)$ and $ b(t,z,P,\varepsilon)$ are analytic functions of
$$|t|\leq T,\ |z|\leq T+1,\ |\varepsilon^KP|\leq
\delta_0=\delta_0(T,K),\ |\varepsilon|\leq\varepsilon_0=\varepsilon_0(T, K)$$
such that $b(t,z,0,\varepsilon)=0$, and
$c(t,z,\varepsilon)$ is an analytic function of
$$|t|\leq T,\ |z|\leq T+1,\ |\varepsilon|\leq\varepsilon_0=\varepsilon_0(T, K).$$
Our goal is to seek a smooth solution $w(t,z)$ of the equation \eqref{yyy}
such that
\begin{center}
$w(t,1)=0$ on $0\leq t\leq T$
\end{center}
 for sufficiently small $\varepsilon$. For completeness, we recall the Nash-Moser theorem as follows.\medskip

\noindent{\textbf{Nash-Moser Theorem}} (see \cite{HM}, p.171,
III.1.1.1) {\textit{Let $\mathfrak{E}_0$ and $\mathfrak{E}$ be tame spaces and 
$\mathfrak{P}:\mathfrak{U}(\subset \mathfrak{E}_0)\to\mathfrak{E}$ a smooth tame map. 
Suppose that the equation for the derivative $D\mathfrak{P}(w)h=g$ has a unique solution 
$h=V\mathfrak{P}(w,g)$ for all $w$ in $\mathfrak{U}$ and all $g$, and that the family of 
inverses $V\mathfrak{P}:\mathfrak{U}\times \mathfrak{E}\to \mathfrak{E}_0$ is smooth tame map. 
Then $\mathfrak{P}$ is locally invertible, and each local inverse $\mathfrak{P}^{-1}$ is 
a smooth tame map.}}\medskip

Now we define the spaces $\mathfrak{E}_0$, $\mathfrak{E}$ and nonlinear mapping $\mathfrak{P}:\mathfrak{U}
(\subset\mathfrak{E}_0)\to
\mathfrak{E}$ by
\begin{align*}
\mathfrak{E}:=&\{ y\in C^{\infty}([-\tau_1,T]\times[0,1]) |\
y(t,z)=0 \  \mbox{for}\ -\tau_1\leq t\leq -\tau_1/2 \},\medskip\\
\mathfrak{E}_0:=&\{ w\in \mathfrak{E}\  |\  w|_{z=1}=0 \},\\
\mathfrak{P}(w):=&\frac{\partial^2 w}{\partial t^2}-\big(1+\varepsilon
a(t,z,-\frac{\partial w}{\partial z},\varepsilon)\big)\triangle w
+\varepsilon b(t,z,-\frac{\partial w}{\partial z}, \varepsilon).
\end{align*}
Here $\tau_1$ is a positive number. Note that $\mathfrak{P}(0)=0$ and take a neighborhood $\mathfrak{U}$ of $0$ such that
$|\varepsilon|^K\|\partial w/\partial z\|_{L^{\infty}}\leq \delta_0 $ for $w \in \mathfrak{U}$.
Then the equation \eqref{yyy} can be written by
\begin{align}\label{owp}
\mathfrak{P}(w)=\varepsilon c(t,z,\varepsilon).
\end{align}
Here we assume $c(t,z,\varepsilon)=0$ for $ t\leq -\tau_1/2$ by change of values of $c$ on $-\tau_1\leq t<0$.
By definition of $\mathfrak{P}$, it is easy to see the Fr\'{e}chet derivative
$D\mathfrak{P}$ of the mapping $\mathfrak{P}$ at a fixed $w\in \mathfrak{U}\subset \mathfrak{E}_0$ is of
the form
$$D\mathfrak{P}(w)h=
\frac{\partial^2 h}{\partial t^2}-\big(1+\varepsilon a_1(t,z,\varepsilon)\big)
\triangle h+
\varepsilon a_2(t,z,\varepsilon)\frac{\partial h}{\partial z}, $$
where
\begin{align}
a_1(t,z,\varepsilon):=&a(t,z,-\frac{\partial w}{\partial z},\varepsilon), \label{a1}\\
a_2(t,z,\varepsilon):=&\frac{\partial a}{\partial P}(t,z,-\frac{\partial w}{\partial z},
\varepsilon)
\triangle w -
\frac{\partial b}{\partial P}(t,z,-\frac{\partial w}{\partial z}, \varepsilon).\label{a2}
\end{align}
Suppose the following statements hold:
\begin{enumerate}
\item[(S1)] $\mathfrak{P}$ is a smooth tame map, $\mathfrak{E}$ being endowed with 
a suitable system of graded norms (for definition, see \cite{HM});
\item[(S2)] for any $w\in \mathfrak{U}\subset \mathfrak{E}_0, g \in \mathfrak{E}$,
there is a unique solution
 $h:=V\mathfrak{P}(w,g)$ of the equation
\begin{align}\label{li}
D\mathfrak{P}(w,h)=g
\end{align}
and the mapping
$V\mathfrak{P}:\mathfrak{U}\times \mathfrak{E}\rightarrow\mathfrak{E}_0$ is a smooth tame map.
\end{enumerate}

\noindent Then it follows from Nash-Moser Theorem that $\mathfrak{P}$ is invertible in a neighborhood $\mathfrak{U}$ of $0$
in $\mathfrak{E}_0$. Thus the inverse image $w =\mathfrak{P}^{-1}(\varepsilon c)$
is a solution of our problem \eqref{yyy} or \eqref{owp}, where $\varepsilon$ is sufficiently small. More precisely, we have the following results.
\newtheorem{theorem}{Theorem}
\begin{theorem}\label{Smoo}
There is a positive constant $\varepsilon_1=\varepsilon_1(T,K)$
such that
for $|\varepsilon|\leq \varepsilon_1$ there exists a smooth
solution $w=w(t,z)$ of \eqref{yyy} defined on
$0\leq t\leq T, 0\leq z\leq 1$ such that $w(t,1)=0$ and
$w=O(\varepsilon)$. In other words,
there is a smooth solution $y=y(t,z)$ of \eqref{laplace} such that $y|_{z=1}=0$
and
$$y(t,z)=y^{(K)}(t,z)+O(\varepsilon^{K+1}). $$
\end{theorem}

\medskip
First of all we must show that the linear equation (51) can be solved uniquely.
But the term $\displaystyle a_2\cdot{\partial h}/{\partial z}$ in
$D\mathfrak{P}(w)h$ could cause trouble, since this term
can have same order as the 
principal part $\displaystyle \triangle h$. When we try to get the energy estimate, keeping in mind that
$(-\triangle h|h)_{\mathfrak{X}}=\|\sqrt{z}\partial h/\partial z\|_{\mathfrak{X}}^2$,
we could not estimate $\|\partial h/\partial z\|_{\mathfrak{X}}$ by
$\|\sqrt{z}\partial h/\partial z\|_{\mathfrak{X}}$ because of the 
singularity at $z=0$. However we have fortunately the following observation:

\begin{proposition}\label{a2z}
For any fixed $w$ in the neighborhood $\mathfrak{U}\subset \mathfrak{E}_0$ there is a smooth function
$\hat{a}_2(t,z,\varepsilon)$ of $0\leq t\leq T, 0\leq z\leq 1, |\varepsilon|
\leq \varepsilon_0(\mathfrak{U})$ such that
$$a_2(t,z,\varepsilon)=z\hat{a}_2(t,z,\varepsilon).$$
\end{proposition}
Proof. By \eqref{a2}, we can write
\begin{equation}
a_2(t,z,\varepsilon)=\frac{\partial a}{\partial P}\cdot z\frac{\partial^2w}{\partial z^2}-
\Big(\frac{N}{2}P\frac{\partial a}{\partial P}+\frac{\partial b}{\partial P}
\Big),
\end{equation}
since
$$\triangle w=z\frac{\partial^2w}{\partial z^2}+
\frac{N}{2}\frac{\partial w}{\partial z}=z\frac{\partial^2w}{\partial z^2}-\frac{N}{2}P.$$
It follows from \eqref{g1}, \eqref{g2} and \eqref{aa}  that
\begin{align*}
\varepsilon\frac{\partial a}{\partial P}=&
\varepsilon^KDG_I(v^{(K)}+\varepsilon^KP)
=\varepsilon^KD^2G(v^{(K)}+\varepsilon^KP),\\
\varepsilon^{K+1}\frac{\partial F_I}{\partial P} =&
\varepsilon^KD^2G(v^{(K)}+\varepsilon^KP)\triangle y^{(K)}\\
=&-\frac{N}{2}D^2G(v^{(K)}+\varepsilon^KP)\cdot v^{(K)}+D^2G(v^{(K)}+\varepsilon^KP)\cdot (z\frac{\partial^2y^{(K)}}{\partial z^2}),\\
\varepsilon^{K+1}\frac{\partial F_{II}}{\partial P}=&
\varepsilon^{K}DG_{II}(v^{(K)}+\varepsilon^KP) \\
=&\varepsilon^K\cdot \frac{N}{2}(v^{(K)}+\varepsilon^KP)\cdot
D^2G(v^{(K)}+\varepsilon^K P),
\end{align*}
since
$$DG_{II}(v)=\frac{N}{2} vD^2G(v),
$$
or
$$\varepsilon\frac{\partial F_{II}}{\partial P}=
\frac{N}{2}(v^{(K)}+\varepsilon^KP)D^2G(v^{(K)}+\varepsilon^KP). $$
Hence, we have
\begin{align*}
\varepsilon\Big(\frac{N}{2}P\frac{\partial a}{\partial P}+\frac{\partial b}{\partial P}
\Big)=&
\varepsilon\Big(\frac{N}{2}P\frac{\partial a}{\partial P}-
\Big(\frac{\partial F_I}{\partial P}
+\frac{\partial F_{II}}{\partial P}\Big)\Big)\\
=&-D^2G(v^{(K)}+\varepsilon^KP)\cdot z\frac{\partial^2y^{(K)}}{\partial z^2}.
\end{align*}
Therefore, $a_2(t,z,\varepsilon)=z\hat{a}_2(t,z,\varepsilon)$ by putting
$$\hat{a}_2(t,z,\varepsilon):=D^2G(v^{(K)}+\varepsilon^KP)
\cdot \frac{\partial^2}{\partial z^2}\varepsilon^{-1}(y^{(K)}+\varepsilon^Kw).$$
$\square$\medskip

Thanks to Proposition \ref{a2z},  we can obtain the following energy inequality.

\begin{proposition} \label{eng}Assume that $|\varepsilon a_1|\leq 1/2$ uniformly for
$t \in \mathbb{R}, \ 0\leq z\leq 1$ and $ |\varepsilon|\leq\varepsilon_0$.
Then there is a constant $C$ such that if 
$h\in \mathfrak{E}_0$ and $ g\in \mathfrak{E}$ satisfy
\begin{equation}\label{62}
\frac{\partial^2 h}{\partial t^2}-(1+\varepsilon a_1)
\triangle h+\varepsilon a_2\frac{\partial h}{\partial z}
=g,
\end{equation}
then
\begin{align}\label{enggg}
||\frac{\partial h}{\partial t}||_{\mathfrak{X}}+
||\sqrt{z}\frac{\partial h}{\partial z}||_{\mathfrak{X}}\leq
C\int_{-\tau_1}^t ||g(\tau)||_{\mathfrak{X}}d\tau\end{align}
for
$0\leq t\leq T$, where
$$||g(t)||_{\mathfrak{X}}=\Big(\int_0^1 |g(t,z)|^2z^{\frac{N}{2}-1}dz
\Big)^{1/2}.$$

\end{proposition}
Proof.  Let us consider the energy
$$E(t):=\int_0^1\Big((h_t)^2+(1+\varepsilon a_1)z(h_z)^2
\Big)z^{\frac{N}{2}-1}dz.$$
We claim that there is a constant $A$ such that
$$E(t)^{1/2}\leq \int_{-\tau_1}^te^{A(t-s)}\|g(s)\|_{\mathfrak{X}}ds.$$
By Proposition \ref{a2z}, the coefficient $a_2(t,z,\varepsilon)$ is of the form
$$a_2(t,z,\varepsilon)=z\hat{a}_2(t,z,\varepsilon),$$
where $\hat{a}_2$ is a smooth function of
$0\leq t\leq T,\ 0\leq z\leq 1$ and $ |\varepsilon|\leq\varepsilon_0$.
Then the equation \eqref{62} turns out to be
\begin{equation}\label{65}
 h_{tt}-(1+\varepsilon a_1)\triangle h+
\varepsilon \hat{a}_2 zh_z=g.
\end{equation}
Multiplying equation \eqref{65} by $h_t$ and integrating it by $d\nu=z^{\frac{N}{2}-1}dz$ from
$z=0$ to $z=1$, we obtain that
\begin{equation}\label{66}
\frac{1}{2}\int^1_0\frac{\partial}{\partial t}h^2_td\nu-\int_0^1(1+\varepsilon a_1)\triangle hh_td\nu+\int_0^1\varepsilon \hat{a}_2zh_zh_td\nu=\int^1_0g(t,z)h_td\nu.
\end{equation}
Using an integration by parts under the boundary condition, we
have
\begin{align}\label{67}
\int_0^1(1+\varepsilon a_1)\triangle hh_td\nu=&\int^1_0(zh_{zz}+\frac{N}{2}h_z)(1+\varepsilon a_1)h_tz^{\frac{N}{2}-1}dz\nonumber\\
=&-\int^1_0h_z\frac{\partial}{\partial z}[(1+\varepsilon a_1)h_tz^{\frac{N}{2}}]dz
+\int_0^1\frac{N}{2}h_z(1+\varepsilon a_1)h_tz^{\frac{N}{2}-1}dz\nonumber\\
=&-\int^1_0h_zh_{zt}(1+\varepsilon a_1)z^{\frac{N}{2}}dz-
\int^1_0\varepsilon h_z (a_1)_zh_tz^{\frac{N}{2}}dz\nonumber\\
=&-\frac{1}{2}\frac{d}{dt}\int^1_0(1+\varepsilon a_1)h_z^2z^{\frac{N}{2}}dz+\frac{1}{2}\int^1_0\varepsilon( a_1)_th_z^2z^{\frac{N}{2}}dz\nonumber\\
&-\int^1_0\varepsilon h_z (a_1)_z h_tz^{\frac{N}{2}}dz.
\end{align}
Since $1+\varepsilon a_1\geq 1/2$, we have $\int_0^1 z(h_z)^2d\nu \leq 2E$. Therefore, it follows from \eqref{65}, \eqref{66} and \eqref{67} that
\begin{align*}
\frac{1}{2}\frac{dE}{dt}=&\frac{1}{2}\varepsilon\int_0^1(a_1)_tz(h_z)^2d\nu
-\varepsilon\int_0^1(a_1)_zzh_zh_td\nu 
-\varepsilon\int_0^1\hat{a}_2zh_zh_td\nu + \int_0^1 gh_td\nu\medskip \\
\leq &AE+ ||g(t)||_{\mathfrak{X}}E^{1/2},
\end{align*}
where
$$A:=\varepsilon(||\partial a_1/\partial t||_{L^{\infty}}+
\sqrt{2}||\sqrt{z}(\partial a_1/\partial z+
\hat{a}_2)||_{L^{\infty}}).
$$
Hence, the Gronwall's argument implies
$$E(t)^{1/2}\leq e^{At}\Big(
E(-\tau_1)^{1/2}+\int_{-\tau_1}^t\|g(s)\|_{\mathfrak{X}}e^{-As}ds
\Big). $$
Since $E(-\tau_1)=0$ from the initial condition for $h\in \mathfrak{E}_0$, 
we get the required inequality \eqref{enggg}. $\square$\medskip

As a corollary of Proposition \ref{eng}, $g=0$ implies $h=0$ by the boundary condition
and which implies that \eqref{li} has a unique solution.
Moreover, this consideration of energy
is sufficient to claim that the inverse $V\mathfrak{P}(w,\cdot)$ of $D\mathfrak{P}(w,\cdot)$ exists.
This can be verified by the standard method on solving the initial boundary value problem
to linear wave equations with smooth coefficients.
See, e.g. Chapter 2 of \cite{Ikawa}.
In fact, for any fixed $t_0$, if we consider the Hilbert space
$\mathfrak{H}=\mathfrak{X}_1\times\mathfrak{X}$ and the operator $\mathfrak{A}(t_0)$, whose domain
$\mathfrak{D}(\mathfrak{A}(t_0))$ is
$$\mathfrak{D}(\mathfrak{A}(t_0))=\{\vec{h}=
(h_0,h_1)^T\in\mathfrak{H} \ | \ h_0\in \mathfrak{X}_2, h_1\in \mathfrak{X}_1, h_0|_{z=1}=h_1|_{z=1}=0\},$$
by $$\mathfrak{A}(t_{0})\begin{pmatrix}
 h_0\\
h_1
\end{pmatrix}=
\begin{pmatrix}
h_1 \\
(1+\varepsilon a_1(t_0,z))\triangle h_0-\varepsilon
\hat{a}_2(t_0,z)h_{0,z}
\end{pmatrix},
$$
then the problem
$$\frac{d\vec{h}}{dt}=\mathfrak{A}(t_0)\vec{h}+\vec{g}(t),
\qquad \vec{h}|_{t=0}=\vec{h}_0\in \mathfrak{D}(\mathfrak{A}(t_0)), $$
where $ \vec{g}(t)=(0,g(t,\cdot))^T$,
allows the application of Hille-Yosida theory. Note that
$-(1+\varepsilon a_1)\triangle h = f$ means that
$$(\sqrt{z}h_z|\sqrt{z}((1+\varepsilon a_1)\phi)_z)_{\mathfrak{X}}=(f|\phi)_{\mathfrak{X}}$$
for any test function $\phi$ or $\phi \in \mathfrak{X}_1$.
Here $\mathfrak{X}_1$ denotes the space of functions $y(z) \in\mathfrak{X}$
such that $\sqrt{z}dy/dz \in \mathfrak{X}$
and $\mathfrak{X}_2$ denotes the space of functions $y \in \mathfrak{X}_1$
such that $-\triangle y\in\mathfrak{X}$. For more details we refer the reader to \cite{Ikawa}.\medskip

Next, we show that the Fr$\acute{\rm e}$chet space $\mathfrak{E}$ is tame for some gradings of norms.
For $y\in \mathfrak{E}$, $n\in\mathbb{N}$, let us define
\begin{align}\label{GN}
&\|y\|_n^{(\infty)}:=
\sup_{0\le j+k\leq n}\Big\|\Big(-\frac{\partial^2}{\partial t^2}
\Big)^{j}(-\triangle)^{k}y\Big\|_{L^{\infty}([-\tau_1,T]\times[0,1])}.
\end{align}
Then we can claim that $\mathfrak{E}$ turns out to be tame by this
grading $(\|\cdot\|_{n}^{(\infty)})_n$ (see \cite{HM}, p.136, II.1.3.6 and p.137, II 1.3.7).
In fact, even if $N$ is not an integer, we can define the Fourier transformation $Fy(\zeta)$ of a function $y(z)$ for $0\leq z<\infty$ by
$$Fy(\zeta):=\int_0^{\infty}K(\zeta z)y(z)z^{\frac{N}{2}-1}dz. $$
Here $K(X)$ is an entire function of $X \in \mathbb{C}$ given by
$$K(X)=2(\sqrt{X})^{-\frac{N}{2}+1}J_{\frac{N}{2}-1}(4\sqrt{X})
=2^{\frac{N}{2}}\Phi_{\frac{N}{2}-1}(X), $$
$J_{\nu}$ being the Bessel function. Then we have
$$F({-\triangle y})(\zeta)=4\zeta\cdot F{y}(\zeta) $$
and the inverse of the transformation $F$ is $F$ itself, see, e.g. \cite{Sneddon}. 
Then it is easy to see $\mathfrak{E}$ endowed with
the grading $(\|y\|_n^{(\infty)})_n$ of the form \eqref{GN}
is a tame direct summand of the tame space
$$L_1^{\infty}(\mathbb{R}\times [0,\infty), d\tau\otimes \zeta^{\frac{N}{2}-1}d\zeta, \log (1+\tau^2+4\zeta)) $$
through the Fourier transformation
$$\mathcal{F}y(\tau, \zeta)=
\frac{1}{\sqrt{2\pi}}\int e^{-\sqrt{-1}\tau t}Fy(t,\cdot)(\zeta)dt$$
and its inverse
applied to the space $C_0^{\infty}((-2T-2\tau_1,2T)\times [0, 2))$, into which functions of
$\mathfrak{E}$ can be extended (see, e.g. \cite{MZ}, p.189, Theorem 3.13) and the space
$$\dot{C}^{\infty}(\mathbb{R}\times [0,\infty)) :=
\{y |\ \forall j, \forall k \lim_{L\rightarrow\infty}\sup_{|t|\geq L,x\geq L}|(-\partial_t^2)^j(-\triangle)^ky|=0\},$$
for which functions of $\mathfrak{E}$ are restrictions.
For the details, see the proof of \cite{HM}, p.137, II.1.3.6.Theorem.

\medskip

On the other hand, let us define
$$\|y\|_n^{(2)}:=\Big(
\sum_{0\le j+k\leq n}
\int_{-\tau_1}^T\|\Big(-\frac{\partial^2}{\partial t^2}\Big)^{j}(-\triangle)^{k}y\|_{\mathfrak{X}}^2dt
\Big)^{1/2}. $$
We have
$$\sqrt{\frac{N}{2}}\|y\|_{\mathfrak{X}}\leq \|y\|_{L^{\infty}}\leq C
\sup_{j\leq\sigma}\|(-\triangle)^jy\|_{\mathfrak{X}},$$
by the Sobolev imbedding theorem (see Appendix A), provided that $2\sigma > N/2$. The derivatives with respect to $t$
can be treated more simply. Then
we see that
the grading $(\|\cdot\|_n^{(2)})_n$ is equivalent to
the grading $(\|\cdot\|_n^{(\infty)})_n$. Hence
$\mathfrak{E}$ is tame with respect to $(\|\cdot\|_n^{(2)})_n$.
The grading $(\|\cdot\|_n^{(2)})_n$ is suitable for energy estimates.
Note that $\mathfrak{E}_0$ is a closed subspace of $\mathfrak{E}$ endowed with these
gradings.
\medskip

 Now we show the statement (S1) by verifying the nonlinear mapping $\mathfrak{P}$
is tame for the grading
$(\|\cdot\|_n^{(\infty)})_n$. To do so,
we write
$$\mathfrak{P}(w)=F(t, z, Dw, w_{tt}, \triangle w), $$
where $D=\partial/\partial z$, $F$ is a smooth function of
$t,z, Dw, w_{tt}, \triangle w$ and linear in $w_{tt}, \triangle w$.
According to \cite{HM} (see p.142, II.2.1.6 and p.145, II.2.2.6), it is sufficient to prove the
linear differential operator $w \mapsto Dw=\partial w/\partial z$ is tame.
But it is clear because of the following result.
\begin{proposition}
For any $m\in\mathbb{N}$ we have the formula
$$\triangle^mDy(z)=
z^{-\frac{N}{2}-m-1}\int_0^z
\triangle^{m+1}y(\zeta)\zeta^{\frac{N}{2}+m}d\zeta. $$
As a corollary it holds that, for any $m,k \in \mathbb{N}$,
$$\|(-\triangle)^mD^ky\|_{L^{\infty}}\leq
\frac{1}{\prod_{j=0}^{k-1}(\frac{N}{2}+m+j)}\|(-\triangle)^{m+k}y\|_{L^{\infty}}.$$
\end{proposition}
Proof. It is easy by integration by parts in induction on $m$ starting from the formula
$$Dy(z)=
z^{-\frac{N}{2}}\int_0^z
\triangle y(\zeta)\zeta^{\frac{N}{2}-1}d\zeta. $$
$\square$\medskip

In parallel with the results of \cite{HM} (see p.144, II.2.2.3.Corollary 
and p.145, II.2.2.5.Theorem), we
should use the following two propositions. Proofs for these propositions are given in Appendix B.

\begin{proposition}\label{p64}
For any positive integer $m$, there is a constant $C$ such that
$$|\triangle^m(f\cdot g)|_0\leq
C(|\triangle^mf|_0|g|_0+
|f|_0|\triangle^mg|_0), $$
where $|\cdot |_0$ stands for $\|\cdot \|_{L^{\infty}}$.
\end{proposition}

\begin{proposition}\label{p65}
Let $F(z,y)$ be a smooth function of $z$ and $y$ and $C_0$ be a positive number. Then
for any positive integer $m$, there is a constant $C>0$ such that
$$|\triangle^mF(z,y(z))|_0\leq C
(1+|y|_m) $$
provided that $|y|_0\leq C_0$, where we denote
$$|y|_m=\sup_{0\leq j\leq m}\|(-\triangle)^jy\|_{L^{\infty}}.$$
\end{proposition}

A tame estimate of the inverse $D\mathfrak{P}(w)^{-1}:g\mapsto h$ will be discussed in the next section. 
This will completes the proof of the main result.

\section{Tame estimate of solutions of linear wave equations}

We consider the wave equation
\begin{equation}
\frac{\partial^2h}{\partial t^2}+\mathcal{A}h=g(t,x), \quad (0\leq t\leq T,\ 0\leq x\leq 1),
\end{equation}
where
\begin{align*}
\mathcal{A}h:=-b_2\triangle h +b_1\check{D}h+b_0h,\quad
\triangle :=x\frac{d^2}{dx^2}+\frac{N}{2}\frac{d}{dx}, \quad \check{D}=x\frac{d}{dx}.
\end{align*}
We denote $\vec{b}(t,x):=(b_2(t,x),b_1(t,x),b_0(t,x))$. The given function $\vec{b}(t,x)$ is supposed to be in
$C^{\infty}([0,T]\times[0,1])$ and we assume that
$|b_2(t,x)-1|\leq 1/2$. The function $g(t,x)$
belongs to $C^{\infty}([0,T]\times[0,1])$ and
we suppose that
\begin{equation}
g(t,x)=0 \ \ \mbox{for}\  0\leq t\leq \tau_2,
\end{equation}
Let us consider the
initial boundary value problem:
\begin{align*}
{\rm  (IBP)}\qquad \left\{
\begin{array}{ll}
\displaystyle\frac{\partial^2h}{\partial t^2}+\mathcal{A}h=g(t,x),\medskip \\
h|_{x=1}=0, \quad
h|_{t=0}=\displaystyle\frac{\partial h}{\partial t}\Big|_{t=0}=0.
\end{array}
\right.
\end{align*}
Then (IBP) admits a unique solution $h(t,x)$ thanks to the energy estimate, and
$h(t,x)=0$ for $0\leq t\leq \tau_2$ because of the uniqueness.
Moreover, since the compatibility conditions are satisfied, the unique
solution turns out to be smooth. A proof can be found e.g. in
\cite[Chapter 2]{Ikawa}.
We are going to get estimates of the higher
derivatives of $h$ by those of $g$ and the coefficients $b_2,b_1,b_0$.\medskip

First, we introduce the following notations:
\begin{notation}\label{nt}~
\begin{enumerate}
\item[{\rm(1)}] For $m,n\in\mathbb{N}$ and for functions $y=y(x)$ of
$x \in [0,1]$, we put
\begin{align*}
(y)_{2m} :=&\|\triangle^my\|, 
\quad \|y\|:=\|y\|_{\mathfrak{X}}:=\Big(\int_0^1|y(x)|^2x^{\frac{N}{2}-1}dx\Big)^{1/2}, \\
(y)_{2m+1} :=&\|\dot{D}\triangle^my\|, \quad \dot{D}=\sqrt{x}\frac{d}{dx}, \\
\|y\|_n:=&\Big(\sum_{0\leq\ell\leq n}(y)_{\ell}^2\Big)^{1/2}, \quad
|y|_n:=\sum_{0\leq\ell\leq n}
\|\dot{D}^{\ell}y\|_{L^{\infty}(0,1)}.
\end{align*}
\item[{\rm(2)}] For $n\in\mathbb{N}$, a fixed $T>0$, and for functions
$y=y(t,x)$ of $(t,x)\in[0,T]\times[0,1]$, we put
\begin{align*}
\|y\|_n^T&:=\Big(\sum_{j+k\leq n}
\int_0^T\|\partial_t^jy\|_k^2dt\Big)^{1/2}, \quad
|y|_n^T:=\sup_{j+k\leq n}\|\partial_t^j\dot{D}^ky\|_{L^{\infty}([0,T]\times[0,1])}.
\end{align*}
Here $\partial_t=\partial/\partial t$. 
\item[{\rm(3)}] Let us say that a grading of norms $(p_n)_{n\in\mathbb{N}}$ is {\bf interpolation admissible} if
for $\ell\leq m\leq n$ it holds
$$p_m(f)\leq Cp_n(f)^{\frac{m-\ell}{n-\ell}}p_{\ell}(f)^{\frac{n-m}{n-\ell}}.
$$
\end{enumerate}
\end{notation}
It is well-known that, $(p_n)_{n\in\mathbb{N}}$ is interpolation admissible if and only if
$$p_n(f)^2\leq C p_{n+1}(f)p_{n-1}(f),\ \ \mbox{
for any }n\geq 1.$$ 
Moreover, if $(p_n)_n$ and $(q_n)_n$ are interpolation admissible, and if
$(i,j)$ lies on the line segment joining
$(k,\ell)$ and $(m,n)$, then
$$p_i(f)q_j(g)\leq C
(p_k(f)q_{\ell}(g)+p_m(f)q_n(g)).
$$ (For a proof, see \cite{HM}, p.144, 2.2.2. Corollary.)\medskip

Since $\dot{D}=\partial/\partial\xi$, where
$x=\xi^2/4$, it is well-known that $(|\cdot|_n)_n$ and
$(|\cdot|_n^T)_n$ are
interpolation admissible. Moreover $(\|\cdot\|_n)_n$ and $(\|\cdot\|_n^T)_n$ are interpolation admissible. To verify it, it is
sufficient to note that
$$y=\sum_{n=1}^{\infty}c_k\phi_k \in C_0^{\infty}[0,1)$$
enjoys
$$(y)_{\ell}=\Big(\sum_k\lambda_k^{\ell}|c_k|^2\Big)^{1/2}. $$
here {$\{\lambda_k\}_{k\in\mathbb{N}}$} are eigenvalues of $-\triangle$ with the Dirichlet
boundary condition at $x=1$ and {$\{\phi_k\}_{k\in\mathbb{N}}$} are associated eigenfunctions.
 Then it is clear by the Schwartz inequality that
$$(y)_n^2\leq (y)_{n+1}(y)_{n-1}$$
for $y \in C_0^{\infty}[0,1)$. Since
$(y)_j\leq (y)_{j'}$ for $j\leq j', \ y\in C_0^{\infty}[0,1)$, we have
$$(y)_{\ell}\leq \|y\|_{\ell}\leq C\cdot (y)_{\ell}\quad\mbox{and}\quad
\|y\|_n^2\leq C\|y\|_{n-1}\|y\|_{n+1}$$
at least for $y \in C_0^{\infty}[0,1)$.
By using a continuous linear extension of functions 
on $[0,1]$ to functions on $[0,2]$ with supports in $[0,3/2)$, we can
claim that this inequality holds for any $y \in C^{\infty}([0,1])$
with a suitable change of the constant $C$. We refer to
\cite{MZ}, Chapter 3, Section 4, Theorem 3.11. It is sufficient to note
the following inequality.
\begin{proposition}
If $\alpha(x) \in C^{\infty}(\mathbb{R})$ is fixed,
then 
$$\|\alpha y\|_n\leq C \|y\|_n.$$
\end{proposition}

A proof can be found in Appendix B. Hence $(\|\cdot\|_n)_n$
and $(\|\cdot\|_n^T)_n$ are interpolation admissible.\\

\subsection{Estimates of the higher
derivatives of $h$}

Our goal is to obtain the following estimates of the higher
derivatives of $h$.

\begin{lemma}\label{hde}
 Assume that $|b_2-1|\leq 1/2$, 
$|\vec{b}|_6^T\leq M$
and $\|g\|_1^T \leq M$. Then there is a constant $C_n=C_n(T,M,N)$ such that
if $h$ is the solution of {\rm (IBP)} then 
$$\|h\|_{n+2}^T\leq C_n
(1+\|g\|_{n+1}^T+|\vec{b}|_{n+5}^T).$$
\end{lemma}

We see that 
$\|y\|_{2m}^T$ is equivalent to
$$\|y\|_m':=\Big(\sum_{j+k\leq m}
\int_0^T(\partial_t^{2j}y)_{2k}^2dt\Big)^{1/2}\quad
\mbox{for }y \in C_0^{\infty}((0,T)\times[0,1)).$$ In fact,  
it is sufficient to note the following inequality.

\begin{proposition}
We have
$$\|\dot{D}\triangle^my\|_{\mathfrak{X}}\leq C(\|\triangle^my\|_{\mathfrak{X}}+
\|\triangle^{m+1}y\|_{\mathfrak{X}}).$$
\end{proposition}

A proof can be found in Appendix B. Therefore $\|\cdot\|_{2m}^T$
is equivalent to $\|\cdot\|_m^{(2)}$. Since $\vec{b}$ is $C^{\infty}$-function of $w, Dw, D^2w$, by Sobolev imbedding theorem, we have
$$ |\vec b|_{n+5}^T\leq C (1+|w|_{n+9}^T)
\leq C'(1+ \|w\|_{n+9+2\sigma}^T), $$
provided $2\sigma>(N+1)/2$. Let $2m=n+2$, then Lemma \ref{hde} reads:
$$ \|h\|_{m}^{(2)}\leq C (1+\|g\|_{n+1}^T+\|w\|_{2m+7+2\sigma}^T) 
              \leq C'(1+\|g\|_m^{(2)}+\|w\|_{m+4+\sigma}^{(2)}). $$ 
Let us sketch a proof of Lemma \ref{hde} in the sequel. 

\subsubsection{Elliptic a priori estimates}

By tedious calculations, we have
\begin{eqnarray*}
[\triangle^m, \mathcal{A}]y :=\triangle^m\mathcal{A}y-\mathcal{A}\triangle^my
=\sum_{j+k=m}(b_{1k}^{(m)}\check{D}\triangle^jy+b_{0k}^{(m)}\triangle^jy),
\end{eqnarray*}
where $\check{D}=xd/dx$ and
\begin{align*}
b_{10}^{(m)}=&-2mDb_2, \\
b_{00}^{(m)}=&-m((2m-1)\triangle +(m-1)(1-N)D)b_2+
m(1+2\check{D})b_1,
\end{align*}
 where $D=d/dx$ and $b_{1k}^{(m)}, b_{0k}^{(m)}, k\geq 1$ are determined by starting
from
\begin{align*}
b_{11}^{(1)}={2Db_0}+(\triangle -(N-2)D)b_1,\quad
b_{01}^{(1)}=\triangle b_0,
\end{align*}
and the {recurrence} formula
\begin{align*}
b_{1k}^{(m+1)}=&b_{1k}^{(m)}+(\triangle-(N-2)D){b_{1k-1}^{(m)}}
+2D{b_{0k-1}^{(m)}}\quad\mbox{for}\quad k\geq 2, \\
b_{11}^{(m+1)}=&b_{11}^{(m)}-4m^2(\triangle+\frac{3-N}{2}D){D}b_2 + \\
&((4m+1)\triangle -(2mN-6m+N-2)D)b_1+2Db_0, \\
b_{0k}^{(m+1)}=& b_{0k}^{(m)}+(1+2\check{D}){b_{1k}^{(m)}}+\triangle {b_{0k-1}^{(m)}}\quad\mbox{for}\quad k\geq 2, \\
b_{01}^{(m+1)}=& b_{01}^{(m)}-m\triangle((2m-1)\triangle + (m-1)(1-N)D)b_2+\\
&m(3+2\check{D})\triangle b_1+\triangle b_0+(1+2\check{D})b_{11}^{(m)}.
\end{align*}
We have used the following calculus formula:
\begin{align*}
\check{D}^2=&x\triangle -\Big(\frac{N}{2}-1\Big)\check{D}, \quad
\triangle\check{D}-\check{D}\triangle =\triangle, \\
\triangle (Q\check{D}P)=&Q\check{D}\triangle P+(1+2\check{D})Q\cdot\triangle P+
(\triangle-(N-2)D)Q\cdot\check{D}P, \\
\triangle(QP)=&Q\triangle P+2(DQ)\check{D}P+(\triangle Q)P.
\end{align*}
It is easy to see that
$$\|b_{0k}^{(m)}\|_{L^{\infty}}\le C|\vec b|_{2k+4}\quad\mbox{and}\quad \|b_{1k}^{(m)}\|_{L^{\infty}}\le  C|\vec{b}|_{2k+2} $$
and therefore
$$\|[\triangle^m, \mathcal{A}]y\|\leq
C\sum_{j+k=m}|\vec{b}|_{2k+4}\|y\|_{2j+1}. $$
Differentiating $[\triangle^m, \mathcal{A}]y$, we get
$$\dot{D}[\triangle^m, \mathcal{A}]y=
\sum_{k+j=m}(\dot{b}_{2k}^{(m)}\triangle^{j+1}y+
\dot{b}_{1k}^{(m)}\dot{D}\triangle^jy+
\dot{b}_{0k}^{(m)}\triangle^jy), $$
where
\begin{align*}
\dot{b}_{2k}^{(m)}=\sqrt{x}b_{1k}^{(m)},\quad
\dot{b}_{1k}^{(m)}=(-\frac{N}{2}+1+\check{D})b_{1k}^{(m)}+b_{0k}^{(m)},\quad
\dot{b}_{0k}^{(m)}=\dot{D}b_{0k}^{(m)}.
\end{align*}
Using
$$\dot{D}\triangle^m[\triangle, \mathcal{A}]=\dot{D}[\triangle^{m+1}, \mathcal{A}]-
\dot{D}[\triangle^m, \mathcal{A}]\triangle, $$
we have
$$\|\dot{D}\triangle^m[\triangle, \mathcal{A}]y\|\leq
C\sum_{j+k=m+1}|\vec{b}|_{2k+6}\|y\|_{2j+2}.$$
Thus 
{
\begin{eqnarray*}
([\triangle,\mathcal{A}]y)_{2m}
\leq C A_m, \quad
([\triangle, \mathcal{A}]y)_{2m+1}\leq C A_m^{\sharp},
\end{eqnarray*}}
where
\begin{align*}
A_m:=\sum_{j+k=m+1}|\vec{b}|_{2k+4}\|y\|_{2j+1},\quad
A_m^{\sharp}:=\sum_{j+k=m+1}|\vec{b}|_{2k+6}\|y\|_{2j+2}.
\end{align*}
Since $A_{m-1}\leq A_{m-1}^{\sharp}\leq A_m\leq A_{m}^{\sharp}$,
we can claim that
{
\begin{align*}
\|[\triangle, \mathcal{A}]y\|_{2m} \leq CA_m,\quad 
\|[\triangle, \mathcal{A}]y\|_{2m+1}\leq CA_m^{\sharp}.
\end{align*}}
 We shall show
\begin{align*}
E(m):\|y\|_{2m+2}\leq& C
(\|\mathcal{A}y\|_{2m}+\|y\|_{2m+1}+
\sum_{j+k=m}|\vec{b}|_{2k+6}\|y\|_{2j+1}), \\
E^{\sharp}(m):\|y\|_{2m+3}\leq&
C(\|\mathcal{A}y\|_{2m+1}+
\|y\|_{2m+2}+
\sum_{j+k=m}|\vec{b}|_{2k+6}\|y\|_{2j+2})
\end{align*}
In fact $E(0)$ is valid, since
$$\triangle y=-\frac{1}{b_2}(\mathcal{A}y-b_1\check{D}y-b_0y) $$
implies 
$$\|\triangle y\|\leq C(\|\mathcal{A}y\|+\|y\|_1). $$
On the other hand,
{
\begin{align*}
\dot{D}\triangle y =&-\frac{1}{b_2}\big(\dot{D}\mathcal{A}y+(-\dot{D}b_2+\sqrt{x}b_1)\triangle y +  (-\frac{N}{2}+1+\check{D}b_1+b_0)\dot{D}y+(\dot{D}b_0)y\big)
\end{align*}}
implies
$$\|\dot{D}\triangle y\|\leq C(\|\dot{D}\mathcal{A}y\|+\|y\|_2)$$
and $E^{\sharp}(0)$ holds. Then $E(m)$ and $E^{\sharp}(m)$ can be verified 
inductively by using the estimates of $\|[\triangle, \mathcal{A}]\|_n$ already obtained.\medskip

By interpolation, we have 
$$\|y\|_{2m+2}\leq C(\|\mathcal{A}y\|_{2m}+\|y\|_{2m+1}+|\vec{b}|_{2m+6}\|y\|_1)$$
and
$$\|y\|_{2m+3}\leq C(\|\mathcal{A}y\|_{2m+1}+\|y\|_{2m+2}+
|\vec{b}|_{2m+6}\|y\|_2),
$$
provided that
$|\vec{b}\|_6 \leq M$. Note that
$$\|y\|_2\leq C(\|\mathcal{A}y\|+\|y\|_1)$$
can be inserted in the estimate of $\|y\|_{2m+3}$.
Thus we have the following result. 
\begin{proposition}\label{prop12}
Suppose $|b_2-1|\leq 1/2$ and $|\vec{b}|_6\leq M$. Then
$$\|y\|_{n+2}
\leq C(\|\mathcal{A}y\|_n+(1+|\vec{b}|_{n+4})\|y\|_1).$$
\end{proposition}

\subsubsection{Estimates for evolutions}

Hereafter we denote generally by
$H$ a solution of the boundary
value problem
$$\frac{\partial^2H}{\partial t^2}+\mathcal{A}H=G(t,x),\qquad H|_{x=1}=0 $$
such that $H(t,x)=0$ for $0\leq t\leq \tau_2$. Thus $\partial_t^jH|_{t=0}=0$ for
any $j\in\mathbb{N}$ and the time derivative $H_j=\partial_t^jH$ satisfies
$$\frac{\partial^2H_j}{\partial t^2} +\mathcal{A}H_j=G_j,\quad H_j|_{x=1}=0,$$
where
$$G_j:=\partial_t^jG-[\partial_t^j,\mathcal{A}]H.$$
We put $G_0=G$.
Provided that $|\vec{b}|_1^T\leq M$, we have the energy
estimate
$$\|\partial_tH\|+\|H\|_1\leq C\int_0^t\|G(t')\|dt'\quad 
\mbox{for }0\leq t\leq T.$$
\begin{remark}
In this subsection, $H$ and $G$ do not mean the
particular functions.
\end{remark}

Let us set
$$Z_n(H):=\sum_{j+k=n}\|\partial_t^jH\|_k.$$
We have the following estimations.
\begin{proposition}\label{estimate}
For $n\in\mathbb{N}$, we have
\begin{enumerate}
\item[{\rm(1)}] 
$Z_{n+2}(H)
\leq C(Z_{n+1}(\partial_tH)+\|G\|_n+
(1+|\vec{b}|_{n+4})\|H\|_1).$
\item[{\rm(2)}] $Z_{n+2}(H)\leq C
(\|\partial_t^{n+1}H\|_1+
\sum\limits_{j+k=n}\|G_j\|_k+\sum\limits_{j+k=n}(1+|\vec{b}|_{k+4})\|\partial_t^jH\|_1).
$
\item[{\rm(3)}] $Z_{n+2}(H)(t)\leq
C\Big(\big(\displaystyle\int_0^tF(t')^2dt'\big)^{1/2}+F(t)\Big), $
where
\begin{align*}
F(t):=&\int_0^t\|\partial_t^{n+1}G(t')\|dt'+
|\vec{b}|_{n+1}^T\|H\|_2^T+\sum_{j+k=n}\|G_j\|_k+\\
&\sum_{j+k=n}(1+|\vec{b}|_{k+4})\|\partial_t^jH\|_1.
\end{align*}
\end{enumerate}
\end{proposition}
{Proof.} (1) By definition, we have
$$Z_{n+2}(H)=Z_{n+1}(\partial_tH)+\|H\|_{n+2}.$$
Then {Proposition \ref{prop12}} implies that 
\begin{align*}
\|H\|_{n+2}&\leq C(\|\mathcal{A}H\|_n+(1+|\vec{b}|_{n+4})\|H\|_1) \leq C(\|\partial_t^2H-G\|_n+(1+
|\vec{b}|_{n+4})\|H\|_1) \\
&\leq C(\|\partial_t^2H\|_n+\|G\|_n+(1+
|\vec{b}|_{n+4})\|H\|_1).
\end{align*}
Note that $\|\partial_t^2H\|_{n}\leq Z_{n+1}(\partial_tH)$. Hence the assertion of this part follows.\medskip

(2) The assertion follows obviously by induction. \medskip

(3) Applying the energy estimate on $\|\partial_t^{n+1}H\|_1$, we get
$$\|\partial_t^{n+1}H\|_1\leq C
\Big(\int_0^t\|\partial_t^{n+1}G\|+\int_0^t
\|[\partial_t^{n+1},\mathcal{A}]H\| \Big).
$$
By interpolation, we have
\begin{align*}
\int_0^t\|[\partial_t^{n+1},\mathcal{A}]H\| &
\leq C\sum_{\alpha+\beta=n+1, \alpha\not=0,}|\partial_t^{\alpha}\vec{b}|_0^t
\Big(\int_0^t
\|\partial_t^{\beta}H\|_2^2\Big)^{1/2} \\
&\leq C'\Big(
\int_0^tZ_{n+2}(H)^2\Big)^{1/2}+|\vec{b}|_{n+1}^T\|H\|_2^T\Big).
\end{align*}
Hence the estimate follows by applying the Gronwall's lemma to
$$Z_{n+2}(H)(t)\leq
C\Big(\big(\int_0^t
Z_{n+2}(H)^2\big)^{1/2}
+F(t)\Big).\ \square$$ 
%
%
%
%

According to Proposition \ref{estimate}, we have
\begin{align}\label{est1}
\|H\|_{n+2}^T\leq & C\Big(\|G\|_{n+1}^T+
\sum_{j+k=n}\Big(\int_0^T\|G_j\|_k^2\Big)^{1/2}+ \nonumber \\
&\sum_{0\leq j\leq n}\sup_{0\leq t\leq T}\|\partial_t^jH\|_1+
|\vec{b}|_{n+1}^T\|H\|_2^T+
|\vec{b}|_{n+4}^T\|H\|_1^T\Big).
\end{align}
Now we estimate the terms in the right-hand side of \eqref{est1}.

\begin{proposition}\label{est2} For $n\in\mathbb{N}$, we have
%
%
\begin{align*}
(1) & \sum_{j+k=n}\Big(\int_0^T\|G_j\|_k^2\Big)^{1/2}\leq C(\|G\|_n^T +\|H\|_{n+1}^T + +|\vec{b}|_{n+4}^T\|H\|_2^T+
|\vec{b}|_{n+5}^T\|H\|_1^T);\\
(2) &
\sup_{0\leq t\leq T}\|\partial_t^jH\|_1 \leq
C(\|G\|_j^T + \|H\|_{j+1}^T 
+|\vec{b}|_{j+4}^T\|H\|_2^T+|\vec{b}|_{j+5}^T\|H\|_1^T).
\end{align*}
\end{proposition}
{ Proof.} (1) It is sufficient to estimate
$$\int_0^T\|[\partial_t^j, \mathcal{A}(\vec{b})]H\|^2dt.$$
Since 
$$[\partial_t^j, \mathcal{A}(\vec{b})]H=
\sum_{\alpha+\beta=j,\alpha\not=0}
\binom{j}{\alpha}
\mathcal{A}(\partial_t^{\alpha}\vec{b})
\partial_t^{\beta}H, $$
and
$$\|\mathcal{A}(\vec{b})y\|_k\leq C
(\|y\|_{k+2}+|\vec{b}|_{k+4}\|y\|_2),$$
we see that 
$$\|\mathcal{A}(\partial_t^{\alpha}\vec{b})\partial_t^{\beta}H\|_k
\leq C(\|\partial_t^{\beta}H\|_{k+2}+
|\partial_t^{\alpha}\vec{b}|_{k+4}\|\partial_t^{\beta}H\|_2,$$
provided that $|\vec{b}|_6^T\leq M$.
By interpolation, we have, for $\alpha+\beta+k=n, \alpha\not=0$,
\begin{align*}
&\Big(\int_0^T\|\mathcal{A}(\partial_t^{\alpha}\vec{b})
\partial_t^{\beta}H\|_k^2\Big)^{1/2}\\
\leq &
C(\|H\|_{\beta+k+2}^T+|\vec{b}|_{\alpha+k+4}^T\|H\|_{\beta+2}^T) 
\\
\leq & C''(\|H\|_{n+1}^T+
|\vec{b}|_5^T\|H\|_{n+1}^T+
|\vec{b}|_{n+5}^T\|H\|_1^T+
|\vec{b}|_{n+4}^T\|H\|_2^T).
\end{align*}
Hence the estimate of this part follows.\medskip

(2) By the energy estimate, we have
$$\|\partial_t^jH\|_1\leq C\int_0^T\|G_j\|.$$
Here we can use the estimate of
$$\int_0^T\|[\partial_t^j, \mathcal{A}]H\|^2$$
given in the proof of the preceding part with
$k=0, n=j$ to obtain the desire estimate. 
$\square$\medskip

\noindent{\bf Proof of Lemma \ref{hde}.} \medskip

Summing up the estimates of Proposition \ref{est2}, we have
\begin{equation}
\|h\|_{n+2}^T\leq
C(\|g\|_{n+1}^T+\|h\|_{n+1}^T+
|\vec{b}|_{n+4}^T\|h\|_2^T +
|\vec{b}\|_{n+5}^T\|h\|_1^T).
\end{equation}
We note that 
$$\|h\|_1^T\leq C\|g\|^T\leq CM.$$
Since 
$$Z_2(h)(t)\leq C
(\|g\|_1^T+\int_0^t
\|[\partial_t, \mathcal{A}]h\|),$$
and
$$\|[\partial_t, \mathcal{A}]h\|=\|\mathcal{A}(\partial_t\vec{b})h\|\leq CZ_2(h).$$
By the Gronwall's lemma, we have
\begin{center}
$Z_2(h)\leq C\|g\|_1^T$\ \  and\ \  $\|h\|_2^T\leq C \|g\|_1^T\leq CM$.
\end{center}
Therefore
we have
\begin{equation}
\|h\|_{n+2}^T\leq C
(\|h\|_{n+1}^T+\|g\|_{n+1}^T+|\vec{b}|_{n+5}),
\end{equation} 
which implies inductively that
\begin{equation}
\|h\|_{n+2}^T\leq C
(1+\|g\|_{n+1}^T+|\vec{b}|_{n+5}).
\end{equation} 
This completes the proof of Lemma \ref{hde}. $\square$.
\\

\noindent {\bf \Large Appendix }\medskip

\noindent{\bf  A. The Sobolev imbedding theorem}\medskip

For the sake of self-containedness, we prove the Sobolev imbedding theorem
for our framework. (The statement is well-known if $N$
is an integer.) Let $y \in \mathfrak{X}$ and $m \in \mathbb{N}, m\geq 1$, we denote
$$((y))_m:=\|(-\triangle)^my\|_{\mathfrak{X}}.$$
Suppose $y\in C_0^{\infty}(0,1)$, then we have the expansion
$$y(z)=\sum_{n=1}^{\infty}c_n\phi_n, $$
where
{$\{\phi_n\}_{n\in\mathbb{N}}$} is the orthonormal system of eigenfunctions of
the operator $T=-\triangle$ with the Dirichlet boundary condition
at $x=1$. Then, for $m\in\mathbb{N}$, we have
$$(-\triangle)^my(x)=\sum_{n=1}^{\infty}c_n\lambda_n^m\phi_n(x)$$
and
$$((y))_m=\Big(\sum_n
|c_n|^2\lambda_n^{2m}
\Big)^{1/2}. $$
As for prerequisites, some properties are illustrated in the sequel.\medskip

\noindent{\bf Lemma A.1.}
Let $j_{\nu,n}$ be the $n$-th positive zero of the Bessel
function $J_{\nu}$, where $\nu=\frac{N}{2}-1$. Then we have
$$\lambda_n=(j_{\nu,n}/2)^2 \sim \frac{\pi^2}{4}n^2 \ \ \mbox{as }n\rightarrow\infty.$$
Proof.
By the Hankel's asymptotic expansion (see \cite{W}), the zeros of $J_{\nu}$ can be determined by
the relation
$$\tan\Big(r-(\frac{\nu}{2}+\frac{1}{4})\pi\Big)=\frac{2}{\nu^2-\frac{1}{4}}r(1+O(r^{-2})).$$
Then we see
$$j_{\nu,n}=\Big(n_0+n+\frac{\nu}{2}+\frac{3}{4}\Big)\pi +
O\Big(\frac{1}{n}\Big) \ \ \mbox{as }n\rightarrow\infty,$$
for some $n_0\in\mathbb{Z}$. $\square$\medskip

\noindent{\bf Lemma A.2.}  There is a constant $C=C(N)$ such that
$$|\phi_n(x)|\leq C n^{\frac{N-1}{2}} \ \ \mbox{for }0\leq x\leq 1.$$
Proof.
We can assume
that $\phi_n(x)$ is the normalization of $\Phi_{\nu}(\lambda_n x)$,
where
$$\Phi_{\nu}\Big(\frac{r^2}{4}\Big)=J_{\nu}(r)\Big(\frac{r}{2}\Big)^{-\nu}.$$
Since $|\Phi_{\nu}(x)|\leq C$ for $0\leq x<\infty$, it is sufficient to estimate
$\|\Phi_{\nu}(\lambda_nx)\|_{\mathfrak{X}}$. Using
the Hankel's asymptotic
expansion in the form
\begin{align*}
J_\nu(r)=&\sqrt{\frac{2}{\pi r}}\Big(\cos\Big(r-\frac{\nu}{2}\pi-\frac{\pi}{4}\Big)
(1+O\Big(\frac{1}{r^2}\Big)\Big)+\\
-&\frac{1}{r}\sin\Big(r-\frac{\nu}{2}\pi-\frac{\pi}{4}\Big)
\Big(\frac{\nu^2-\frac{1}{4}}{2}+O\Big(\frac{1}{r^2}\Big)\Big)\Big),
\end{align*}
 we see that
\begin{align*}
\|\Phi_{\nu}(\lambda_nx)\|_{\mathfrak{X}}^2=&
(\lambda_n)^{-\nu-1}\int_0^{j_{\nu,n}}J_{\nu}(r)^2rdr
= (\lambda_n)^{-\nu-1}\Big(\frac{1}{\pi}j_{\nu,n}+O(1)\Big) \\
=&(\lambda_n)^{-\nu-1}\cdot \frac{2}{\pi}(\lambda_n^{1/2}+O(1))
\sim
\frac{2}{\pi}(\lambda_n)^{-\nu-\frac{1}{2}}.
\end{align*}
Then Lemma A.1 implies that
$$\|\Phi_{\nu}(\lambda_nx)\|_{\mathfrak{X}}^{-1}\sim \mbox{Const.}\cdot n^{\nu+\frac{1}{2}}.$$
$\square$\medskip

\noindent{\bf Lemma A.3.} If $ y \in C_0^{\infty}(0,1)$ and $0\leq j\leq m$, then
$((y))_j\leq((y))_m.$\medskip

\noindent Proof.
For
$y=\sum c_n\phi_n$, we have
\begin{align*}
((y))_j^2=&
\sum |c_n|^2\lambda_n^{2j} =(\lambda_1)^{2j}\sum |c_n|^2(\lambda_n/\lambda_1)^{2j} \\
\leq& (\lambda_1)^{2j}\sum|c_n|^2(\lambda_n/\lambda_1)^{2m}=
\lambda_1^{2j-2m}((y))_m^2.
\end{align*}
According to \cite{W} (see Section 15-6 , p.208), we know that $j_{\nu,1}$ is an increasing function of
$\nu>0$ and $j_{\frac{1}{2}, 1}=\pi$. Therefore,
$\lambda_1\geq (\pi/2)^2 >1$ for $N\geq 2$ and which implies $((y))_j\leq((y))_m$. $\square$\medskip

\noindent{\bf Lemma A.4.} If $2s>N/2$, then
there is a constant $C
=C(s,N)$ such that
$$\|y\|_{L^{\infty}}\leq C((y))_s\quad\mbox{
for any }y\in C_0^{\infty}(0,1).$$
\noindent Proof.
Let $y=\sum c_n\phi_n(x)$, then Lemmas A.1 and A.2 imply
that
\begin{align*}
|y(x)|\leq
\sum |c_n||\phi_n(x)|
\leq
C\sum |c_n|n^{\frac{N-1}{2}} 
\leq C\sqrt{\sum|c_n|^2\lambda_n^{2s}}
\sqrt{\sum n^{N-4s-1}}.
\end{align*}
Since $N-4s<0$, the last term in the above inequality is finite. Therefore we get the required estimate. 
$\square$\medskip

Now, for $R>0$, we denote by $\mathfrak{X}(0,R)$ the Hilbert space of functions $y(x)$
of $0\leq x\leq R$ endowed with the inner product
$$(y_1|y_2)_{\mathfrak{X}(0,R)}=
\int_0^Ry_1(x)\overline{y_2(x)}x^{\frac{N}{2}-1}dx. $$
Moreover, for $m\in\mathbb{N}$, we denote by $\mathfrak{X}^{2m}(0,R)$ the space of functions $y(x)$ of
$0\leq x\leq R$ for which the derivatives
$(-\triangle)^jy \in \mathfrak{X}$ exist in the sense of distribution for $0\leq j\leq m$. And we use the norm
$$\|y\|_{\mathfrak{X}^{2m}(0,R)}:=\Big(
\sum_{0\leq j\leq m}\|(-\triangle)^jy\|_{\mathfrak{X}(0,R)}^2
\Big)^{1/2}. $$
Let us denote by $\mathfrak{X}_0^{2m}(0,R)$ the closure of $C_0^{\infty}(0,R)$ in the space
$\mathfrak{X}^{2m}(0,R)$. There is a continuous linear extension $\Psi: \mathfrak{X}^{2m}(0,1) \rightarrow
\mathfrak{X}_0^{2m}(0,2)$ such that
$$\|y\|_{\mathfrak{X}^{2m}(0,1)}\leq \|\Psi y\|_{\mathfrak{X}^{2m}(0,2)}
\leq C\|y\|_{\mathfrak{X}^{2m}(0,1)}.$$
See \cite{MZ}, p.186, Theorem 3.11, keeping in mind
Propositions 6,7. Then,
 by Lemmas A.3 and A.4, the Sobolev imbedding theorem holds for
$y \in \mathfrak{X}_0^{2s}(0,2)$. Say, if $2s >N/2$, there is a constant $C$
such that $$\|y\|_{L^{\infty}}\leq C\|y\|_{\mathfrak{X}^{2s}(0,2)} $$
for $y \in \mathfrak{X}_0^{2s}(0,2)$. Thus the same imbedding theorem holds for
$y \in C^{\infty}([0,1]) \subset \mathfrak{X}^{2s}(0,1)$
through the above extension. The conclusion is that, if $2s >N/2$, there is a constant
$C=C(s,N)$ such that

$$\|y\|_{L^{\infty}}\leq C\sup_{0\leq j \leq s}\|(-\triangle)^jy
\|_{\mathfrak{X}} $$
for any $y \in C^{\infty}([0,1])$.\\

\noindent{\bf  B. Nirenberg-Moser type inequalities}\medskip

Let us prove Propositions 8--11
in the sequel.\\

\noindent{\bf Proof of Proposition 8.}\medskip

First, it is easy to verify the formula
$$\dot{D}^kDy(z)=z^{-\frac{N+k}{2}}\int_0^z\dot{D}^k\triangle y(z')
(z')^{\frac{N+k}{2}-1}dz',\eqno{\rm(B.1)}
$$
where $k\in \mathbb{N}$,
$$\dot{D}:=\sqrt{z}\frac{d}{dz}\quad\mbox{and} \quad D:=\frac{d}{dz}. $$
Since $\triangle =\dot{D}^2+\frac{N-1}{2}D$, (B.1) implies
$$|\dot{D}^kDy|_0\leq \frac{2}{N+k}|\dot{D}^{k+2}y|_0+\frac{N-1}{N+k}|\dot{D}^kDy|_0. $$
Here and hereafter $|\cdot|_0$ stands for $\|\cdot\|_{L^{\infty}}$. Thus we have
$$|\dot{D}^kDy|_0\leq \frac{2}{k+1}|\dot{D}^{k+2}y|_0.$$
Repeating this estimate, we get
$$|\dot{D}^kD^jy|_0\leq \Big(\frac{2}{k+1}\Big)^j|\dot{D}^{k+2j}y|_0. \eqno{\rm(B.2)}$$
On the other hand, since
$\dot{D}^2=\triangle -\frac{N-1}{2}D$ and $D\triangle -\triangle D=D^2$, we have
$$\dot{D}^{2\mu}=\sum_{k=0}^{\mu}C_{k\mu}\triangle^{\mu-k}D^k \eqno{\rm(B.3)}$$
with some constants $C_{k\mu}=C(k,\mu,N)$.
Then it follows from (B.3) and Proposition 7 that
$$|\dot{D}^{2\mu}D^jy|_0\leq C|\triangle^{\mu+j}y|_0. \eqno{\rm(B.4)}$$
Since
$$\triangle =\dot{D}^2+\frac{N-1}{2}D\quad\mbox{and}\quad D\dot{D}^2-\dot{D}^2D=D^2, $$
it is easy to see that there are constants
$C_{km}=C(k,m,N)$ such that
$$\triangle^m=\sum_{k=0}^m
C_{km}\dot{D}^{2(m-k)}D^k. \eqno{\rm(B.5)}$$
Applying the Leibnitz' rule to $D$ and $\dot{D}$, we see
$$\triangle^m(f\cdot g)=
\sum C_{k\ell jm}(\dot{D}^{2(m-k)-\ell}D^{k-j}f)\cdot(\dot{D}^{\ell}D^jg) \eqno{\rm(B.6)}$$
with some constants $C_{k\ell jm}$. The summation is taken for
$0\leq j\leq k\leq m, 0\leq \ell \leq 2(m-k)$. By estimating
each term of the right-hand side of
(B.6), we can obtain the assertion of Proposition 8. In fact, we consider
the term
$$(\dot{D}^{\ell'}D^{j'}f)\cdot(\dot{D}^{\ell}D^jg)$$
provided that
$\ell'+\ell +2(j'+j)=2m$. By (B.2) and (B.4) we have
\begin{eqnarray*}
|\dot{D}^{\ell}D^jg|_0\leq
C|\dot{D}^{\ell +2j}g|_0
\leq  C'|\dot{D}^{2m}g_0|^{\frac{\ell +2j}{2m}}|g|_0^{1-\frac{\ell+2j}{2m}}
\leq  C''|\triangle^mg|_0^{\frac{\ell+2j}{2m}}|g|_0^{1-\frac{\ell+2j}{2m}}
\end{eqnarray*}
for some positive constants $C$, $C'$ and $C''$.
Here we have used the Nirenberg interpolation for $\dot{D}=\partial/\partial \xi$,
where $x=\xi^2/4$. The same
estimate holds for
$|\dot{D}^{\ell'}D^{j'}f|_0$. Therefore we have
\begin{align*}
|(\dot{D}^{\ell'}D^{j'}f)\cdot(\dot{D}^{\ell}D^jg)|_0 \leq &
 C|\triangle^mf|_0^{\frac{\ell'+2j'}{2m}}|f|_0^{1-\frac{\ell'+2j'}{2m}}|\triangle^mg|_0^{\frac{\ell+2j}{2m}}
|g|_0^{1-\frac{\ell+2j}{2m}} \\
\leq & C(|\triangle^mf|_0|g|_0+|f|_0|\triangle^mg|_0),
\end{align*}
since $X^{\theta}Y^{1-\theta}\leq X+Y$. $\square$\\

\noindent{\bf Proof of Proposition 9.}\medskip

Suppose $F(z,y)$ is a smooth function of $z$ and $y$. Let us consider
the composed function $U(z):=F(z,y(z))$. We claim that
$$|\triangle^mU|_0\leq C(1+|y|_m) $$
provided that $|y|_0\leq C_0$.
 In fact,
$$\triangle^mU=\sum C_{km}\dot{D}^{2(m-k)}D^k U $$
consists of several terms of the following form:
$$\Big(\dot{D}_x^K\Big(\frac{\partial}{\partial y}\Big)^LD_x^k\Big(\frac{\partial}{\partial y}\Big)^{\ell}F\Big)\cdot
(\dot{D}^{K_1})\cdots(\dot{D}^{K_L}y)\cdot(\dot{D}^{\mu_1}D^{k_1}y)\cdots
(\dot{D}^{\mu_{\ell}}D^{k_{\ell}}y), $$
where
$$k+k_1+\cdots+k_{\ell}=\kappa,$$ $$ K+K_1+\cdots K_L+\mu_1+\cdots+\mu_{\ell}=2(m-\kappa).$$
Therefore
$$K_1+\cdots+K_L+(\mu_1+2k_1)+\cdots(\mu_{\ell}
+2k_{\ell}) \leq 2m. $$
Applying the Nirenberg interpolation to $\dot{D}$ and using (B.4),
we have
$$|\dot{D}^{K_1}y|_0\leq C|y|_m^{\frac{K_1}{2m}}|y|_0^{1-\frac{K_1}{2m}}. $$
Similarly,
$$|\dot{D}^{\mu_1}
D^{k_1}y|_0\leq C|\dot{D}^{\mu_1+2k_1}y|_0\leq C'|y|_m^{\frac{\mu_1+2k_1}{2m}}
|y|_0^{1-\frac{\mu_1+2k_1}{2m}}, $$
and so on. Then our claim follows obviously. $\square$\\

We note that by (B.2), (B.4) and (B.5) we have
$$\frac{1}{C}|\dot{D}^{2j}f|_0\leq |\triangle^jf|_0
\leq C|\dot{D}^{2j}f|_0.$$

\noindent{\bf Proof of Proposition 10.}\medskip

It can be verified that
$$\triangle^m(\alpha y)=
\sum_{j+k=m}(\alpha_{1k}^{(m)}\check{D}\triangle^jy+
\alpha_{0k}^{(m)}\triangle^jy), \eqno{\rm (B.7)}$$
where
$\alpha_{1k}^{(m)}$ and $\alpha_{0k}^{(m)}$ are determined by the {recurrence} formula
\begin{align*}
\alpha_{1k}^{(m+1)}=&\alpha_{1k}^{(m)}+
(\triangle -(N-2)D){\alpha_{1k-1}^{(m)}}+2D{\alpha_{0k-1}^{(m)}},\\
\alpha_{0k}^{(m+1)}=&(1+2\check{D})\alpha_{1k}^{(m)}+
\alpha_{0k}^{(m)}+\triangle\alpha_{0,k-1}^{(m)},
\end{align*}
starting from
$$\alpha_{10}^{(0)}=0, \qquad \alpha_{00}^{(0)}=\alpha. $$
Here we have used the convention $\alpha_{1k}^{(m)}=\alpha_{0k}^{(m)}
=0$ for $k<0$ or $k>m$. Of course $\alpha_{10}^{(m)}=0$ for any $m$.
Therefore we see that {$\|\triangle^m(\alpha y)\|\leq C \|y\|_{2m}$.}
Differentiating the formula (B.7), we get
$$
\dot{D}\triangle^m(\alpha y) =
\sum_{j+k=m}(\dot{\alpha}_{2k}^{(m)}\triangle^{j+1}y+\dot{\alpha}_{1k}^{(m)}\dot{D}\triangle^jy +
\dot{\alpha}_{0k}^{(m)}\triangle^jy),
$$
where
\begin{eqnarray*}
\quad\dot{\alpha}_{2k}^{(m)}=\sqrt{x}\alpha_{1k}^{(m)},\quad
\dot{\alpha}_{1k}^{(m)}=\Big(-\frac{N}{2}+1+\check{D}\Big)\alpha_{1k}^{(m)}+
\alpha_{0k}^{(m)}, \quad
\alpha_{0k}^{(m)}=\dot{D}\alpha_{0k}^{(m)}.
\end{eqnarray*}
It is clear that {$\|\dot{D}\triangle^m(\alpha y)\|
\leq C\|y\|_{2m+1}$}, since $\dot{\alpha}_{20}^{(m)}=0$
for any $m$. $\square$\\

\noindent{\bf Proof of Proposition 11.}\medskip

It is sufficient to prove that
$$\|\dot{D}y\|\leq C(\|y\|+\|\triangle y\|), $$
where and hereafter we denote $\|\cdot\|=\|\cdot\|_{\mathfrak{X}}$.
If $w$ satisfies the Dirichlet boundary condition $w(1)=0$, then
$$\|\dot{D}w\|^2=(-\triangle w\ |\ w)\leq \|\triangle w\|\|w\|.$$
Therefore
we have
$$\|\dot{D}y\|^2\leq \|\triangle y\|(\|y\|+|y(1)|).$$
On the other hand we have
$$\sqrt{\frac{2}{N}}|y(1)|\leq \|y\|+\sqrt{\frac{2}{N-2}}\|\dot{D}y\|.$$
In fact, since
$$y(1)=y(z)+\int_z^1\frac{1}{\sqrt{z'}}\dot{D}y(z')dz',$$
we have
$$|y(1)|^2\leq |y(z)|^2 +\frac{2}{N-2}\|\dot{D}y\|^2z^{-\frac{N}{2}+1}$$
for $z>0$. Integrating this, we get the above estimate of $|y(1)|$.
Hence we have, for any $\epsilon>0$,
\begin{align*}
\|\dot{D}y\|^2&\leq C\|\triangle y\|(\|y\|+\|\dot{D}y\|) 
\leq C\Big(\frac{1}{2\epsilon}\|\triangle y\|^2+\frac{\epsilon}{2}(\|y\|+\|\dot{D}y\|)^2\Big) \\
&\leq C\big(\frac{1}{2\epsilon}\|\triangle y\|^2+\epsilon \|y\|^2+
\epsilon \|\dot{D}y\|^2\Big).
\end{align*}
Taking $\epsilon$ to be small, we get the desired estimate. $\square$


\begin{thebibliography}{3}
\bibitem{CS} D. Coutand and S. Shkoller, Well-posedness in smooth function spaces for
moving-boundary 1-d compressible Euler equations in physical vacuum, {\it Comm. Pure Appl. Math.} 
{\bf LXIV} (2011), pp. 328-366.

\bibitem{DS} N. Dunford and J. T. Schwartz, {\it Linear Operators}, Part II, Wiley, 1963, NY.

\bibitem{HM} R. Hamilton, The inverse function theorem of Nash and Moser, {\it Bull. American Math. Soc.} {\bf 7} (1982), pp. 65-222.

\bibitem{HLM1} C.-H. Hsu, S.-S. Lin and T. Makino, Periodic solutions to the 1-dimensional compressible Euler equation with gravity, {\it Hyperbolic Problems-
Theory, Numerics and Applications}, Yokohama Publishers, 2006, pp. 163-170.

\bibitem{HLM2} C.-H. Hsu, S.-S. Lin and T. Makino, Smooth solutions to a class of quasilinear wave equations, 
{\it J. Diff. Eqns.} {\bf 224} (2006), pp. 229-257.

\bibitem{Ikawa} M. Ikawa, {\it  Hyperbolic Partial Differential Equations and
Wave Phenomena} (Translations of Math. Monographs, Vol. 189), AMS, Providence, Rhode Island, 2000.

\bibitem{JM} J. Jang and N. Masmoudi, Well-posedness for compressible Euler equations
with physical vacuum singularity, {\it Comm. Pure Appl. Math.} {\bf LXII} (2009), pp. 1327-1385.

\bibitem{Liu} T.-P. Liu, Compressible flow with damping and vacuum,
{\it Japan J. Appl. Math.} {\bf 13} (1996), pp. 25-32.

\bibitem{LiuY} T.-P. Liu and T. Yang, Compressible flow with vacuum and physical singularity,
 {\it Methods Appl. Anal.} {\bf 31} 
(2000), pp. 223-237.


\bibitem{MZ} S. Mizohata, {\it The Theory of Partial Differential Equations},
Cambridge University Press, 1973.

\bibitem{Sneddon} I. N. Sneddon, {\it Fourier Transforms}, NY, McGraw-Hill, 1951; NY, Dover, 1995.


\bibitem{W}  G. N. Watson, {\it A Treatise on the Theory of Bessel Functions}, Cambridge University Press, 1958.

\bibitem{Y} T. Yang, Singular behavior of vacuum states
for compressible fluids, {\it
Comput. Appl. Math.}
{\bf 190} (2006), pp. 211-231.
\end{thebibliography}
\end{document}